\documentclass[11pt,leqno]{article}

\usepackage[all]{xy}				
\usepackage{amsthm}				
\usepackage{amsfonts}				
\usepackage{bussproofs}				
\usepackage{cancel}				
\usepackage[shortlabels]{enumitem}
\usepackage{fancyhdr}				
\usepackage[stable]{footmisc}			
\usepackage{geometry}				
\usepackage{graphicx}				
\usepackage{hyperref}				
\usepackage[utf8]{inputenc}			
\usepackage{lineno}				
\usepackage{mathrsfs}				
\usepackage{MnSymbol}				
\usepackage[numbers]{natbib}			
\usepackage{pdfpages}				
\usepackage{tikz}				
\usetikzlibrary{tikzmark}
\usepackage{url}				
\usepackage{xcolor}				

\theoremstyle{definition}       		
\newtheorem{df}{Definition}	    		
\newtheorem{lm}[df]{Lemma}	  	    	
\newtheorem{prop}[df]{Proposition}		
\newtheorem{rmrk}[df]{Remark}			
\newtheorem{tr}[df]{Theorem}	 		

	\newif\ifSuppressEndOfProva\SuppressEndOfProvafalse

\newcommand{\noi}{\noindent}
\newcommand{\setl}{\setlength\itemsep{-0.2em}}

\newcommand{\oid}{\ensuremath{\Leftrightline}}
\newcommand{\noid}{\ensuremath{\nLeftrightline}}

\newcommand{\C}{\ensuremath{\mathcal{C}}}
\newcommand{\E}{\ensuremath{\textbf{E}}}
\newcommand{\D}{\ensuremath{\textbf{D}}}
\newcommand{\K}{\ensuremath{\mathcal{K}}}

\newcommand{\sig}{\ensuremath{\mathcal{L}}}
\newcommand{\ov}{\ensuremath{\overline}}
\newcommand{\wh}{\ensuremath{\widehat}}

\newcommand{\fde}{$F\!D\!E$}

\newcommand{\nf}{$N\!\mathit{4}$}
\newcommand{\nt}{$N\!\mathit{4}$}

\newcommand{\qfde}{$Q\!F\!D\!E$}
\newcommand{\qvfde}{$Q_{v}F\!D\!E$}
\newcommand{\qvletf}{$Q_{v}LET_{F}^-$}
\newcommand{\qnf}{$Q\!N\!\mathit{4}$}
\newcommand{\nm}{$N^{-}$}
\newcommand{\ffde}{$F\!F\!D\!E$}
\newcommand{\fnf}{$F\!N\!\mathit{4}$}

\newcommand{\red}{\color{red}}

\newcommand{\az}{\color{blue}}

\begin{document}


\title{\textbf{On Universally Free First-Order Extensions of Belnap-Dunn's
Four-Valued Logic and Nelson's Paraconsistent Logic \nf}\footnote{This and a
forthcoming paper from the authors on the origins and motivations of the
interpretation of \fde\ and \nf\ in terms of inconsistent and incomplete
information states are complementary. Here the technical details are given in
full, while historical and conceptual discussions are kept to a minimum. In the
second paper, on the other hand, the technical results are summarized and the
philosophical and historical aspects are discussed in much more detail. We would
like thank Luiz Carlos Pereira for comments that helped to shape some of the
ideas in this paper.}}

\author{
	Henrique Antunes\footnote{Corresponding author.}\\
	Federal University of Bahia\\
	Department of Philosophy\\
	Salvador, Brazil\\
	\texttt{henrique.antunes@ufba.br}
	\\\and
	Abilio Rodrigues\footnote{The research of the second author is supported
	by the grants 310037/2021-2, 408040/2021-1, and 402705/2022-0, Conselho
	Nacional de Desenvolvimento Científico e Tecnológico (CNPq, Brazil), and
	APQ-02093-21, Fundação de Amparo à Pesquisa do Estado de Minas Gerais
	(FAPEMIG, Brazil).}\\
	Federal University of Minas Gerais\\
	Department of Philosophy\\
	Belo Horizonte, Brazil\\
	\texttt{abilio.rodrigues@gmail.com}
}

\date{}

\maketitle


	\begin{abstract}

	The aim of this paper is to introduce the logics \ffde\ and \fnf, which
	are universally free versions of Belnap-Dunn's four-valued logic, also
	known as the logic of first-degree entailment (\fde), and Nelson's
	paraconsistent logic \nm\ (a.k.a. \qnf). Both \fde\ and \qnf\ are
	suitable to be interpreted as information-based logics, that is, logics
	that are capable of representing the deductive behavior of possibly
	inconsistent and incomplete information in a database. Like \qnf\ and
	some non-free first-order extensions of \fde, \ffde\ and \fnf\ are
	endowed with Kripke-style variable domain semantics, which allows
	representing the dynamic aspect of information processing, that is, how
	a database receives new information over time, including information
	about new individuals. We argue, however, that \ffde\ and \fnf\ can
	better represent the development of inconsistent and incomplete
	information states (i.e., configurations of a database) over time than
	their non-free versions. First, because they allow for empty domains,
	which corresponds to the idea that a database may acknowledge no
	individual at all at an early stage of its development. Second, because
	they allow for empty names, which get interpreted as information about
	new individuals is inserted into the database. Also, both systems
	include an identity predicate that is interpreted along the same lines
	of the other logical operators, viz., in terms of independent positive
	and negative rules.

	\end{abstract}

\noi \textbf{Keywords:} Belnap-Dunn's four-valued logic; first-degree
entailment; Nelson's paraconsistent logic \nf; information-based logics; free logics

\section{Introduction}\label{sec:introduction}

In the final section of \cite{Belnap1977}, the 1977 paper that introduces
the four-valued semantics for \textit{the logic of first-degree entailment}
(\fde), alongside with a corresponding intuitive interpretation in terms of
incomplete and inconsistent information, Belnap writes:

	\begin{quote}

	Quantifiers introduce a number of subtleties to which I shall merely tip
	my hat, while recognizing that treating them in detail is quite
	essential to my enterprise$\dots$In any event, {the semantics given for
	the connectives extend to universal and existential quantifiers in an
	obvious way, and I suppose the job done}. And the various alternatives
	mentioned above turn out not to make any difference to the
	\textit{logic}. (p.~50)

	\end{quote}

\noi Belnap is right in recognizing the importance of introducing quantifiers to
\fde. After all, the task of formulating a logic intended to express the
deductive behavior of inconsistent and incomplete information could hardly be
considered to be fully accomplished without such an extension. However, it is
fair to say that Belnap underestimated the challenges posed by this task, and
that his conclusion that the job was done in the passage above was in fact too
hasty. In effect, there are different ways of introducing quantifiers to \fde\
\cite{Antunes&Rodrigues2024,Antunesetal2022,Beall&Logan2017,Priest2002}, and it is far
from obvious how to do so while preserving Belnap's original ideas concerning
\fde.

As an attempt to meet the challenge posed by Belnap's comments, in this paper we
present a universally free version of \fde, dubbed here \ffde, together with a
corresponding adequate variable domain Kripke semantics. We will argue that
\ffde\ enjoys several properties that make it a good candidate for an
information-based first-order extension of \fde. Since the task of extending
\fde\ could hardly be considered complete without a suitable treatment of
identity, in addition to $\forall$ and $\exists$, \ffde\ also includes the
identity predicate \oid, whose account differs significantly from other accounts
given elsewhere in the literature \cite{Antunesetal2022,Beall&Logan2017}. In
spite of having several desirable properties, \ffde\ has, however, a drawback
that is inherited from \fde: it lacks an implication connective that validates
\textit{modus ponens}. So, in the final part of the paper we will also present
the variable domain free version with identity \fnf\ of \textit{Nelson's
paraconsistent logic} \nf, which improves upon \ffde\ by extending it with a
constructive implication.

The method used here is that of conducting a technical investigation of both
systems while taking into account their intended intuitive interpretations, as
well as some relevant historical aspects of their development. Our guiding idea
is to combine the technical results with the interpretation of \fde\ as the
underlying logic for reasoning in scenarios with inconsistent and incomplete
information. We will argue that allowing empty domains and languages with
uninterpreted constants, together with variable domains, is better suited for
the intended intuitive interpretations of \fde\ and \nf\ than the semantics
formulated so far. Specifically, these aspects allow one to represent a sequence
of configurations of a database -- or \textit{information states} -- starting
with a possibly empty database and a first-order language with any number of
constants, some of which may be uninterpreted. Then, as information about new
individuals is inserted into the database, they get assigned to the previously
empty constants. As for the non-standard treatment of identity proposed here, we
want identity formulas and their negations to be independent from each other,
allowing for scenarios in which $c_{1} \oid c_{2}$ and $c_{1} \noid c_{2}$ hold,
and for scenarios in which neither do. Thus, we will adopt independent semantic
clauses for $\oid$ and $\noid$, as well as independent positive and negative
inference rules.

The remainder of this paper is structured as follows. In Section
\ref{sec:desiderata} we go over some of the desiderata that a logic intended to
deal with incomplete and inconsistent information should satisfy,viz., (i) it
must be capable of representing positive and negative information in such a way
that they are independent from each other; (ii) it must be capable of
representing the dynamic aspect of information processing; (iii) it must make
room for the possibility of new individuals being acknowledged at further stages
of the development of a database; (iv) it must allow for empty databases, as
well as for databases whose linguistic resources exceed the individuals it
countenances at a given time. In addition to arguing that \ffde\ satisfies
(i)-(iv), we also provide some preliminary explanations of the ingredients that
make up the logic. In Section \ref{sec:ffde} we present both a natural deduction
system (Subsection \ref{sec:natural_deduction}) and a corresponding free
variable domain Kripke semantics for \ffde\ (Subsection \ref{sec:semantics}).
The final part of Section \ref{sec:ffde} (Subsection
\ref{sec:soundness_and_completeness}) contains detailed proofs of soundness and
completeness results for \ffde. In Section \ref{sec:fn4} we then show how \ffde\
can be modified to yield the logic \fnf, which is a free version with identity
of the logic \nm\ of Almukdad and Nelson \cite{Almukdad&Nelson1984}.

\section{On information-based first-order extensions of \fde\footnote{Parts of
sections 2.1 and 2.2 have already appeared in \cite{Antunes&Rodrigues2024} and
\cite{Antunesetal2022}.}}\label{sec:desiderata}

\subsection{\fde\ as an information-based logic}\label{sec:fde}

Belnap's \cite{Belnap1977} motivation for proposing a four-valued semantics for
\fde\ was to design a logic to be used by a computer, which was conceived by him
to be a ``question-answering system'' that answers questions based on deductions
that take what is stored in its database as premises \cite[p.~36]{Belnap1977}.
So interpreted, \fde\ is an \textit{information-based logic}, in the sense of a
logic suitable for processing information, governing which conclusions may be
drawn from the information stored in the computer's database. Since the
databases envisaged by Belnap may contain contradictory information about a
certain topic or lack of information thereof, the computer is supposed to be
able to handle these situations in a sensible way, and so its underlying logic
must be both paraconsistent and paracomplete (see \cite[p.~35]{Belnap1977}).

One of the core ideas of Belnap's interpretation of \fde\ is the distinction
between \textit{positive} and \textit{negative information}. According to this
view, \textit{every} piece of information that is processed by the computer has
a ``polarity'' that indicates whether its content has been affirmed or denied.
Specifically, given a sentence $A$, \textit{positive information} $A$ is the
information that $A$ is true ($A$ has been affirmed), while \textit{negative
information} $A$, or, equivalently, positive information $\neg A$, is the
information that $A$ is false ($A$ has been denied).

Another core idea is that positive and negative information are independent from
each other, in the sense that the presence of positive (negative) information
$A$ does not exclude the presence of negative (positive) information $A$, nor
does the absence of positive (negative) information $A$ ensures the presence of
negative (positive) information $A$. This idea is closely connected to the
previous one. After all, the fact that a certain piece of information $A$ has
been affirmed does not exclude the possibility of its also being denied -- e.g.,
$A$ may have been affirmed by source $S_{1}$ while being denied by source
$S_{2}$. Moreover, the fact that a certain piece of information $A$ has not
been affirmed does not mean that it has been denied, for none of the sources of
the computer may have information about $A$.

The four-valued semantics of \fde, whose values are \textit{\sf T}, \textit{\sf
F}, \textit{\sf B}, and \textit{\sf N}, is able to express these two ideas in a
seamless way. According to Belnap \cite[p.~38]{Belnap1977}, \textit{\sf T} and
\textit{\sf F} are to be taken as `told true' and `told false' signs in the
sense that a computer ``has been told'' that $A$ is true (has received positive
information $A$), and that it has been told that $A$ is false (has received
negative information $A$). The value \textit{\sf B} means that the computer has
been told that $A$ is both true and false (it has received both positive and
negative information $A$), while \textit{\sf N} means that nothing about $A$ has
been told to the computer (it has received no information about $A$). In view of
this interpretation, the following four scenarios emerge:

	\begin{itemize}\setl \label{prop.4.scenarios}

	\item[] $v(A) = \sf T$: $A$ holds and $\neg A$ does not hold: only
	positive information $A$;

	\item[] $v(A) = \sf F$: $\neg A$ holds and $A$ does not hold: only
	negative information $A$;

	\item[] $v(A) = \sf B$: both $A$ and $\neg A$ hold: contradictory
	information;

	\item[] $v(A) = \sf N$: neither $A$ nor $\neg A$ holds: no information
	at all.

	\end{itemize}

A third important idea is that a database may evolve over time, acquiring new
information at further stages of its development.  Although this idea is
suggested by some of the comments made by Belnap in \cite{Belnap1977}, he
did not incorporate it into the semantics of \fde, nor did he take it into
account while describing how to extend \fde\ with quantifiers (see
\cite[p.~36]{Belnap1977}). As has been already argued by Wansing
\cite{Wansing1993}, this dynamic aspect of information processing may be represented
within the framework of a suitable Kripke semantics, and in this paper we shall
follow Wansing's lead while formulating an information-based first-order
extension of \fde.\footnote{Wansing \cite{Wansing1993} extends the informal reading
of Kripke semantics for intuitionistic logic to logics capable of expressing
positive and negative information, conceived of as independent from each other.
Thus, he writes ``information states may develop differently depending of the
basic information acquired in the course of time'' (6) and ``[a]ny theory of
information processing, in order to be viewed as adequate, will be expected to
allow for representing both positive as well as negative information'' (p.~13).}

\subsection{First-order extensions of \fde}\label{sec:qfde}

As mentioned in the Introduction, \textit{our} task of formulating a first-order
extension of \fde\ is not that of merely coming up with a sound and complete
system. Rather, it is to do so while remaining faithful to Belnap's original
proposal of formulating a well-motivated information-based logic. As it turns
out, this task is not as straightforward as Belnap expected.

Quantified versions of \fde\ have already appeared in
\citep{Anderson&Belnap1963,Antunesetal2022,Beall&Logan2017,Priest2002}, all of
which validate the ``inferential version'' of the \textit{constant domain
axiom}:

	\begin{align}
	\tag{CD} \forall x(B \lor A) \vDash B \lor \forall x A
	\end{align}

\noi (where $x$ is not free in $B$). Thus, any corresponding adequate deductive
system must either have (CD) as a primitive inference rule or replace the usual
$\forall$-introduction rule

	\begin{center}

		\bottomAlignProof
				\AxiomC{$A(c/x)$}
			\RightLabel{$\forall I$}
			\UnaryInfC{$\forall xA$}
		\DisplayProof

	\end{center}

\noi by its stronger version

	\begin{center}

	\bottomAlignProof
			\AxiomC{$B \lor A(c/x)$}
		\RightLabel{$\forall I_{c}$}
		\UnaryInfC{$B \lor \forall xA$}
	\DisplayProof \label{rule.intro.CD}

	\end{center}

\noi which mimics the classical multiple conclusion $\forall$-right
rule.\footnote{See Remark 29 of \citep{Antunesetal2022} for a proof of this fact.}

As its name itself implies, (CD) is valid in semantics that take the universal
quantifier to have a constant domain, that is, a domain of individuals that does
not vary across the stages of a model -- either because they are required to
share the same domain, or because the semantics is not even a Kripke semantics
to begin with. On the other hand, (CD) may fail to hold in models whose stages
may have different domains, that is, it is invalid in a variable domain
semantics.

We have already seen that a Kripke semantics is to be preferred over a purely
extensional semantics due to its capacity to represent the dynamic aspect of
information processing, viz., an information state turning into another by
having some new pieces of information added to it. However, if the stages of a
model are required to share the same domain, then no new information can concern
an individual whose existence was previously unacknowledged. After all, since
all stages have the same domain, there can be no ``new'' individuals, but only
new information about the same old individuals. Since we believe that it may be
useful to allow the information that gets added to a database to be about
previously unacknowledged individuais, such as when one is told that a certain
newborn baby has just left the hospital, we think that a Kripke semantics with
variable domains is to be preferred over a single domain semantics.

Now, in \citep{Antunes&Rodrigues2024} we presented a variable domain first-order version of
\fde, dubbed \qvfde. Unlike the first-order extensions of \fde\ in
\citep{Anderson&Belnap1963,Antunesetal2022,Beall&Logan2017,Priest2002}, \qvfde\
does not validate (CD) and its $\forall$-introduction rule is $\forall I$ rather
than $\forall I_{c}$. The logic \ffde\ to be introduced below is a universally
free version with identity of \qvfde, and, concerning the task of formulating an
information-based first-order extension of \fde, it improves upon \qvfde\ in at
least three respects: in addition to having positive and negative rules for
$\oid$, by allowing empty domains it makes it possible to represent stages of
the development of a database at which it concerns no individual at all. Also,
by allowing empty terms, \ffde\ can represent scenarios in which the linguistic
resources of the computer exceed the individuals acknowledged in its database at
a certain stage.

\subsection{Free logics}\label{sec:free_logics}

Free logics are a family of systems that were proposed and investigated from the
1950s on by Leonard, Lambert, Hintikka, van Fraassen, among others (see e.g.
\citep{Hintikka1959,Lambert1960,Leonard1956,vanFraassen&Lambert1967}).\footnote{For
comprenhensive surveys of free logics, see
\cite{Bencivenga2002,Lehmann2002,Nolt2007}.} The name `free logic', coined by
Lambert \citep{Lambert1960}, stands for `logic free of existential
presuppositions in general and with respect to singular terms in particular',
which reveals that these systems have strong ontological motivations. More
specifically, the proponents of free logics subscribed to the traditional
existential reading of the quantifiers, while rejecting that individual terms
invariably denote existent objects (see \cite[p.~148-8]{Bencivenga2002},
\cite[p.~197]{Lehmann2002}, \cite[p.~2]{Morscher&Heike2001}).

From a technical point of view, free logics allow \textit{empty singular terms},
i.e., terms to which the interpretation function (of a model) assigns no element
of the quantificational domain $D$.\footnote{This can be achieved in two
alternative ways: either by assigning individuals that do not belong to $D$ to
empty terms (\textit{dual-domain semantics}) or by assigning them no object at
all (\textit{single domain semantics}).} As a result, the inference rules below
are invalid in the vast majority of free logics:

	\begin{center}

			\bottomAlignProof
					\AxiomC{$\forall xA$}
				\RightLabel{$\forall E$}
				\UnaryInfC{$A(c/x)$}
			\DisplayProof
		\qquad
			\bottomAlignProof
					\AxiomC{$A(c/x)$}
				\RightLabel{$\exists I$}
				\UnaryInfC{$\exists xA$}
			\DisplayProof

	\end{center}

\noi but can be recovered provided that they include an additional premise that
expresses that $c$ is not an empty constant:

	\begin{center}

			\bottomAlignProof
					\AxiomC{$\forall xA$}
					\AxiomC{$\E c$}
				\RightLabel{$\forall E_{\E}$}
				\BinaryInfC{$A(c/x)$}
			\DisplayProof
		\qquad
			\bottomAlignProof
					\AxiomC{$A(c/x)$}
					\AxiomC{$\E c$}
				\RightLabel{$\exists I$}
				\BinaryInfC{$\exists xA$}
			\DisplayProof

	\end{center}

\noi (where \E\ is the \textit{existence predicate}).

In some presentations of free logic, $\E c$ is defined as $\exists x(x \oid c)$,
but it can also be regarded as a primitive logical predicate alongside $\oid$,
in which case $\E c$ and $\exists x(x \oid c)$ turn out to be equivalent
provided that the logic in question is sententially classical, i.e., validates
all tautological inferences. The fact that \E\ is usually interpreted as
expressing existence and the equivalence between $\E c$ and $\exists x(x \oid
c)$ are symptomatic of the aforementioned ontological interpretation of free
logic. It is worth  emphasizing, however, that nothing in the semantics or
proof theory of a free logic enforces such an interpretation: \E\ and $\exists$
may be taken to be ontologically neutral expressions, neither presupposing nor
implying existence. In fact, this is precisely how we will understand these
symbols in the context of \ffde, but in order to avoid potential confusions we
will write \D\ instead of \E\ and call the former `definedness predicate'. The
rationale for this view is that even though $\D c$ may be interpreted as
expressing existence in a given database, in the sense that $c$ denotes an
object that exists \textit{according to the database}, this does not mean that
the putative object assigned to $c$ by the interpretation function
\textit{really} exists.

Unlike other predicates of a first-order logic, \D\ will be taken to behave
classically in \ffde: for a given constant $c$, one and at most one of the
sentences $\D c$ and $\neg \D c$ will hold in a given stage of a model. Hence,
the following rules will be valid in \ffde:

	\begin{center}

			\bottomAlignProof
					\AxiomC{$\D c$}
					\AxiomC{$\neg \D c$}
				\RightLabel{$PEX_{\D}$}
				\BinaryInfC{$B$}
			\DisplayProof
		\qquad
			\bottomAlignProof
					\AxiomC{}
				\RightLabel{$PEM_{\D}$}
				\UnaryInfC{$\D c \lor \neg \D c$}
			\DisplayProof

	\end{center}

\noi This behavior stems from the particular interpretation of $\D c$ and $\neg
\D c$ to be adopted here: $\D c$ holds if the interpretation $I(c,w)$ of $c$ at
stage $w$ is defined, while $\neg \D c$ holds if $I(c,w)$ is undefined.
Moreover, since (the negation of) an atomic sentence with an empty constant
conveys no information (for every such sentence must be about individuals whose
existence is acknowledged), the following rules will be valid as well:

	\begin{center}

			\bottomAlignProof
					\AxiomC{$Pc_{1}\dots c_{m}$}
				\RightLabel{$\D I$}
				\UnaryInfC{$\D c_{i}$}
			\DisplayProof
		\!\!
			\bottomAlignProof
					\AxiomC{$\neg Pc_{1}\dots c_{m}$}
				\UnaryInfC{$\D c_{i}$}
			\DisplayProof

	\end{center}

\noi for every $1 \leq i \leq m$.

A free logic that also for allows models with empty domains is called
\textit{inclusive} or \textit{universally free}. Unlike non-inclusive free
logics, whose $\forall$-introduction and $\exists$-elimination rules are
typically the same as those of classical first-order logic, in an inclusive
logic these rules assume the following format:

	\begin{center}

			\bottomAlignProof
					\AxiomC{$[\E c]^{i}$} \noLine
					\UnaryInfC{$\vdots$} \noLine
					\UnaryInfC{$A(c/x)$}
				\RightLabel{$\forall I_{i}$}
				\UnaryInfC{$\forall xA$}
			\DisplayProof
		\qquad
			\bottomAlignProof
					\AxiomC{$\exists xA$}
						\AxiomC{$[\E c, A(c/x)]^{i}$}
						\UnaryInfC{$\vdots$} \noLine
					\UnaryInfC{$B$}
				\RightLabel{$\exists E_{i}$}
				\BinaryInfC{$B$}
			\DisplayProof

	\end{center}

\noi Since \ffde\ is universally free, its $\forall$-introduction and
$\exists$-elimination rules must undergo similar changes, as so must its
$\neg \exists$-introduction and $\neg \forall$-elimination rules.

\subsection{Identity}\label{sec:identity}

As we have seen above, a core aspect of Belnap's interpretation of \fde\
as an information-based logic is the distinction between positive and negative
information, which are taken to be independent from each other and required to
be neither exclusive nor exhaustive. These ideas manifest themselves both in the
semantics of \fde, which includes independent positive and negative clauses for
each connective, and in the corresponding deductive systems, by including
positive and negative natural deduction rules. A first-order extension of \fde\
that is meant to preserve this core aspect of Belnap's interpretation should
also be expected to have independent semantic clauses and deduction rules for
the quantifiers, and if this extension includes the identity predicate $\oid$,
then the treatment of $\oid$ should follow the same lines.

All the first-order extensions of \fde\ presented in the literature so far,
viz.,
\citep{Anderson&Belnap1963,Antunes&Rodrigues2024,Antunesetal2022,Beall&Logan2017,Priest2002},
have both positive and negative semantic clauses and deduction rules for the
quantifiers, but among these only the systems in \cite{Antunesetal2022} and
\cite{Beall&Logan2017} include the identity predicate. In Beall and Logan's
\cite{Beall&Logan2017} semantics, $\oid$ is a classical predicate, satisfying
both $\vDash c_{1} \oid c_{2} \lor c_{1} \noid c_{2}$ and $c_{1} \oid c_{2},
c_{1} \noid c_{2} \vDash B$ (see \cite[Remark 34]{Antunesetal2022} for further
details). As for the system in \cite{Antunesetal2022}, i.e., \qfde, which is a
constant domain first-order extension of \fde, identity formulas are not
required to satisfy the schemas $\vDash A \lor \neg A$ and $A, \neg A \vDash B$
-- and so \qfde\ allows for models in which both $c_{1} \oid c_{2}$ and $c_{1}
\noid c_{2}$ hold, as well as scenarios in which neither do (see \cite[Remark
26]{Antunesetal2022} for further details). As far as the task of formulating an
information-based first-order extension of \fde\ is concerned, we take \qfde's
treatment of identity to be an improvement over the one in
\cite{Beall&Logan2017}. After all, there is no reason for supposing that the
third and fourth scenarios mentioned in Subsection \ref{sec:fde} do not arise
w.r.t. to information conveyed by identity formulas. However, the treatment of
identity in \cite{Antunesetal2022} does not go far enough, for although \qfde\
has positive rules for identity, it lacks negative rules, telling us nothing
about when one is authorized to infer formulas of the form $c_{1} \noid c_{2}$.
Among other respects, \ffde\ improves upon \qfde\ by addressing this issue.

Ignoring some philosophical issues concerning identity, it is fair to say that
the principles of \textit{indiscernibility of identicals} and of
\textit{identity of indiscernibles}, respectively,

	\begin{enumerate}[label={(\arabic*)}, resume]\setl

	\item If $a$ is identical to $b$, then for every property $\textbf{P}$,
	$a$ has $\textbf{P}$ iff $b$ has $\textbf{P}$,\label{prop.indisc.ident}

	\item If for every property $\textbf{P}$, $a$ has $\textbf{P}$ iff $b$
	has $\textbf{P}$, then $a$ is identical to $b$,\footnote{Henceforth, we
	shall restrict our attention to properties and to unary predicates. It
	is straightforward to generalize the considerations below to $n$-ary
	relations and $n$-ary predicates. Moreover, while talking about the
	semantics of \ffde, we will ignore the fact that such notions as
	extension, anti-extension, holds etc. must be relativized to a stage of
	a model. Thus, we shall write `$a$ belongs to the extension of
	$\textbf{P}$' instead of `$a$ belongs to the extension of $\textbf{P}$
	at $w$', `$\textbf{P}c$ holds' instead of `$\textbf{P}c$ holds at $w$',
	and so and so forth.} \label{prop.ident.indisc}

	\end{enumerate}

\noi provide a conceptual characterization of identity.\footnote{This
characterization is explicitly incorporated in second-order logic by means of
the definition of $t_{1} \oid t_{2}$ as $\forall P(Pt_{1} \leftrightarrow
Pt_{2})$. Note, however, that this definition requires all instances of the
comprehension schema, viz., $\exists P \forall x(Px \leftrightarrow A)$, to be
valid. Hence, the definition won't work in non-standard formulations of
secord-order logic.} The contrapositives of \ref{prop.indisc.ident} and
\ref{prop.ident.indisc},

	\begin{enumerate}[label={(\arabic*)}, resume]\setl

	\item If for some property $\textbf{P}$, $a$ has $\textbf{P}$ but $b$
	does not have $\textbf{P}$ (or vice-versa), then $a$ is not identical to
	$b$, \label{prop.neg.indisc.ident}

	\item If $a$ is not identical to $b$, then for some property
	$\textbf{P}$, $a$ has $\textbf{P}$ but $b$ does not have $\textbf{P}$
	(or vice-versa), \label{prop.neg.ident.indisc}

	\end{enumerate}

\noi are known in the literature as, respectively, the principles of
\textit{diversity of the dissimilar} (see, e.g., \cite{Williams2008}) and of
\textit{dissimilarity of the diverse} (see, e.g., \cite{Forrest2020}).
Clearly, \ref{prop.indisc.ident} and \ref{prop.neg.indisc.ident} are
(classically) equivalent, as are \ref{prop.ident.indisc} and
\ref{prop.neg.ident.indisc}.

If `not' behaves classically, and so `not-$\textbf{P}$' designates the
complement of $\textbf{P}$, then \ref{prop.neg.indisc.ident} is equivalent to:

	\begin{enumerate}[label={(\arabic*)}, resume]\setl

	\item If for some property $\textbf{P}$, $a$ has $\textbf{P}$ and $b$
	has not-$\textbf{P}$ (or vice-versa), then $a$ is not identical to
	$b$.\label{prop.neg.indisc.ident2}

	\end{enumerate}

\noi But the equivalence fails to hold when `not' is paraconsistent and
paracomplete. In effect, if `not' is paraconsistent, then one is not authorized
to infer that $b$ does not have $\textbf{P}$ from the fact that it has
not-$\textbf{P}$, for $b$ may have \textit{both} $\textbf{P}$ and
not-$\textbf{P}$ (and so (3) does not entail (5)); and if `not' is paracomplete,
then one is not authorized to infer that $b$ has not-$\textbf{P}$ from the fact
that it does not have $\textbf{P}$, for $b$ may have neither $\textbf{P}$ nor
not-$\textbf{P}$ (and so (5) does not entail (3)).

An alternative way of expressing \ref{prop.indisc.ident} and
\ref{prop.neg.indisc.ident} that avoids the talk of objects and properties is by
respectively:

	\begin{enumerate}[label={(\arabic*)}, resume]\setl

	\item If $Pc_{1}$ and $c_{1} \oid c_{2}$ hold, then $Pc_{2}$ also holds,
	for every predicate letter $P$;\label{prop.indisc.ident.alternative}

	\item If $Pc_{1}$ and $\neg Pc_{2}$ hold, for some predicate letter $P$,
	then $c_{1} \noid c_{2}$
	holds.\label{prop.neg.indisc.ident.alternative}\footnote{Notice
	that the corresponding versions of \ref{prop.ident.indisc} and
	\ref{prop.neg.ident.indisc}, viz.,

		\begin{itemize}\setl

		\item[] If for every predicate letter $P$, $Pc_{1}$ holds iff
		$Pc_{2}$ holds, then $c_{1} \oid c_{2}$ holds;

		\item[] If $c_{1} \noid c_{2}$, then for some predicate $P$,
		$Pc_{1}$ holds and $Pc_{2}$ does not hold (or vice-versa),

		\end{itemize}

	\noi do not hold in classical first-order logic, for there are models with
	different objects that satisfy exactly the same predicates.}

	\end{enumerate}

\noi As before, \ref{prop.indisc.ident.alternative} and
\ref{prop.neg.indisc.ident.alternative} are equivalent when $\neg$ is classical,
\ref{prop.indisc.ident.alternative} may fail to entail
\ref{prop.neg.indisc.ident.alternative} when $\neg$ is paraconsistent, and
\ref{prop.neg.indisc.ident.alternative} may fail to entail
\ref{prop.indisc.ident.alternative} when $\neg$ is paracomplete. Hence, these
two principles should be expected to be independent from each other in a
first-order extension with identity of \fde. Now, we believe that not only
should \ref{prop.indisc.ident.alternative} hold in \ffde, for the obvious reason
that it represents one of the most basic properties of identity, but so does
\ref{prop.neg.indisc.ident.alternative}, since it lays down very natural
sufficient conditions for one to obtain information that certain objects are not
the same. Fortunately, it is possible to retain both
\ref{prop.indisc.ident.alternative} and \ref{prop.neg.indisc.ident.alternative}
by adopting independent semantic clauses and deduction rules for $\oid$ and
$\noid$.

As in \qfde\ and \qvfde, in \ffde\ predicate letters are interpreted in terms of
\textit{extensions} and \textit{anti-extensions}: given a predicate letter $P$,
its extension $P_{+}$ is the set of individuals of the domain $D$ that satisfy
$Px$, while its anti-extension $P_{-}$ is the set of individuals of $D$ that
satisfy $\neg Px$. $P_{+}$ and $P_{-}$ are independent from each other in the
sense that they are required to be neither exclusive nor exhaustive w.r.t.  $D$
-- that is, $P_{+} \cap P_{-}$ may be different from $\emptyset$ and $P_{+} \cup
P_{-}$ may be a proper subset of $D$. Since $\oid$ is a binary predicate letter,
we should expect it to be interpreted in a similar way, viz., in terms of an
extension $\oid_{+}$ and an anti-extension $\oid_{-}$. However, in order to
ensure that \ref{prop.indisc.ident.alternative} and
\ref{prop.neg.indisc.ident.alternative} hold for \ffde\ we need to impose some
specific conditions on $\oid_{+}$ and $\oid_{-}$. In the case of $\oid_{-}$, it
suffices to require it to include every pair $\langle a, b \rangle$ such that
$a$ belongs to $P_{+}$ and $b$ belongs to $P_{-}$ (or vice-versa), for some
predicate letter $P$. As for $\oid_{+}$, it may seem that the standard
interpretation of identity would give us everything we need.  After all, by
taking $\oid_{+}$ to be the identity relation on $D$, $\oid$ would satisfy
\ref{prop.indisc.ident.alternative}, as well as other standard properties of
identity, e.g., \textit{reflexivity}. But due to some technical issues, which
also arise in connection with the behavior of identity in intuitionistic logic,
we are prevented from taking this route.\footnote{For more about these issues,
see \cite[\S\S 5.11-12]{Troelstra&vanDalen1988} and Remark \ref{rmrk:congruence}
below.} Instead, we will interpret $\oid_{+}$ as a \textit{congruence relation}
w.r.t.  to the extensions and anti-extensions of predicate letters. That is,
rather than being the set $\{\langle a, a \rangle\!: a \in D\}$, $\oid_{+}$ will
be any equivalence relation $\sim$ on $D$ such that $a \sim b$ only if (i) $a
\in P_{+}$ iff $b \in P_{+}$ and (ii) $a \in P_{-}$ iff $b \in P_{-}$, for every
$P$. This requirement ensures that both \ref{prop.indisc.ident.alternative} and
the reflexivity of $\oid$ hold for \ffde.

Positive and negative natural deduction rules for $\oid$ can be easily obtained
from the considerations above. Rule $\oid\!E$ below is one of $\oid$'s positive rules,
which in \ffde\ must be formulated in terms of atomic formulas and their
negations:

	\begin{center}

			\bottomAlignProof
					\AxiomC{$Pc_{1}$}
					\AxiomC{$c_{1} \oid c_{2}$}
				\RightLabel{$\oid\!E$}
				\BinaryInfC{$Pc_{2}$}
		\!\!
			\DisplayProof
			\bottomAlignProof
					\AxiomC{$\neg Pc_{1}$}
					\AxiomC{$c_{1} \oid c_{2}$}
				\BinaryInfC{$\neg Pc_{2}$}
			\DisplayProof

	\end{center}

\noi $\noid\!I$ is $\oid$'s single negative rule:

	\begin{center}

			\bottomAlignProof
					\AxiomC{$Pc_{1}$}
					\AxiomC{$\neg Pc_{2}$}
				\RightLabel{$\noid\!I$}
				\BinaryInfC{$c_{1} \noid c_{2}$}
			\DisplayProof

	\end{center}

\noi As we shall see, both $\oid\!E$ and $\noid\!I$ can be generalized to
non-atomic formulas (see Proposition \ref{prop:igeneralization} below).

Finally, the usual rule of introduction of identity, $\oid$'s second positive
rule, also holds, but it must be restricted to constants that have been defined:

	\begin{center}

			\bottomAlignProof
					\AxiomC{$\D c$}
				\RightLabel{$\oid\!I_{\D}$}
				\UnaryInfC{$c \oid c$}
			\DisplayProof

	\end{center}

	\begin{rmrk} It is worth noting that the contrapositive of $\oid\!I$
	cannot be valid if we are to allow for contradictions involving $\oid$.
	Indeed, $c_{1} \noid c_{1}$ follows from $c_{1} \oid c_{2}$ and $c_{1}
	\noid c_{2}$ by $\oid\!E$. Since $\D c_{1}$ follows from $c_{1} \oid
	c_{2}$ by $\D I$, one cannot be able to infer $\neg \D c_{1}$ from
	$c_{1} \noid c_{1}$, for otherwise any pair of formulas of the form
	$c_{1} \oid c_{2}$ and $c_{1} \noid c_{2}$ would lead to triviality (by
	rule $PEX_{\D}$).

	\end{rmrk}

\section{The logic \ffde}\label{sec:ffde}

In the previous section we discussed the desiderata that a first-order
extension of \fde\ must satisfy in order to deliver a well-motivatived
information-based logic that is faithful to Belnap's original ideas in
\cite{Belnap1977}: (i) it must be capable of representing positive and
negative information, which are taken to be independent from each other and are
required to be neither exclusive nor exhaustive; (ii) it must be capable of
representing the dynamic aspect of information processing, that is, represent
the possibility of a database having new pieces of information inserted into it
over time; (iii) it must make room for the possibility of new individuals being
acknowledged at further stages of the development of a database; (iv) it must
allow for empty databases, as well as for databases whose linguistic resources
exceed the individuals it countenances at a given time.

Now that we have gone over each of these desiderata and given a preliminary
explanation of why \ffde\ is capable of satisfying them, we will present the
logic in much more detail. Specifically, in subsections
\ref{sec:natural_deduction} and \ref{sec:semantics} we will present a natural
deduction system and a corresponding free variable domain Kripke semantics for
\ffde. In Subsection \ref{sec:soundness_and_completeness} we will then prove
that the system in Subsection \ref{sec:natural_deduction} is both sound and
complete.

\subsection{Natural Deduction System}\label{sec:natural_deduction}

The logical vocabulary of \ffde\ is composed by the connectives $\neg$, $\land$,
and $\lor$, the quantifiers $\forall$ and $\exists$, the individual variables
from $\mathcal{V} = \{v_{i}: i \in \mathbb{N} \}$, and parentheses. We shall
take a \textit{first-order language} \sig\ to be a pair $\langle \mathcal{C},
\mathcal{P} \rangle$ such that \C\ is an infinite set of \textit{individual
constants} and $\mathcal{P}$ is a set of \textit{predicate letters}.
$\mathcal{P}$ is not required to be non-empty, and each element $P$ of
$\mathcal{P}$ is assumed to have a finite arity. We will denote the set
$\mathcal{P} \cup \{\oid\}$ by $\mathcal{P}^{\oid}$.

We assume here the usual definitions of such notions as \emph{formula},
\emph{bound/free occurrence of a variable}, \emph{sentence}, etc. Given a
first-order language \sig, the set of formulas and the set of sentences
generated by \sig\ will be denoted by $Form(\sig)$ and $Sent(\sig)$,
respectively. $x$, $x_{1}$, $x_{2}$ will be used as metavariables ranging over
$\mathcal{V}$; $c$, $c_{1}$, $c_{2}$,$\dots$, as metavariables ranging over
$\mathcal{C}$; $A$, $B$, $C$, $\dots$, as metavariables ranging over
$Form(\sig)$; and $\Gamma$, $\Delta$, $\Sigma$, as metavariables ranging over
subsets of $Form(\sig)$. Given $c \in \C$, $A(c/x)$ is the formula that results
by replacing every free occurrence of $x$ in $A$ by $c$.

In order to avoid certain technical difficulties that would arise with open
formulas, the deductive systems presented below are formulated exclusively with
respect to sentences. This is the reason why we have assumed right from the
outset that first-order languages must always have an infinite stock of
individual constants -- for otherwise we could be prevented from applying some
of the quantifier rules due to the lack of enough constants.

	\begin{df}\label{df:nds} Let $\sig = \langle \C, \mathcal{P} \rangle$ be
	a first-order language, $c,c_{1},\dots,c_{m} \in \C$, $P \in
	\mathcal{P}$, and $A, B, C \in Sent(\sig)$. The logic \ffde\ is defined
	over \sig\ by the following natural deduction rules:

		{\footnotesize

		\begin{center}

				\bottomAlignProof
					\AxiomC{$A$}
						\AxiomC{$B$}
					\RightLabel{$\land I$}
					\BinaryInfC{$A \land B$}
				\DisplayProof
			\qquad
				\bottomAlignProof
						\AxiomC{$A \land B$}
					\RightLabel{$\land E$}
					\UnaryInfC{$A$}
				\DisplayProof
			$\!\!$
				\bottomAlignProof
						\AxiomC{$A \land B$}
					\UnaryInfC{$B$}
				\DisplayProof
			\qquad
				\bottomAlignProof
						\AxiomC{$A$}
					\RightLabel{$\lor I$}
					\UnaryInfC{$A \lor B$}
				\DisplayProof
			$\!\!$
				\bottomAlignProof
						\AxiomC{$B$}
					\UnaryInfC{$A \lor B$}
				\DisplayProof
			\qquad
				\bottomAlignProof
						\AxiomC{$A \lor B$}
							\AxiomC{$[A]^{i}$} \noLine
						\UnaryInfC{$\vdots$} \noLine
						\UnaryInfC{$C$}
							\AxiomC{$[B]^{i}$} \noLine
						\UnaryInfC{$\vdots$} \noLine
						\UnaryInfC{$C$}
					\RightLabel{$\lor E^{i}$}
					\TrinaryInfC{$C$}
				\DisplayProof

		\end{center}

		\vspace{0.25cm}

		\begin{center}

				\bottomAlignProof
						\AxiomC{$\neg A$}
					\RightLabel{$\neg \land I$}
					\UnaryInfC{$\neg(A \land B)$}
				\DisplayProof
			$\!\!$
				\bottomAlignProof
						\AxiomC{$\neg B$}
					\UnaryInfC{$\neg(A \land B)$}
				\DisplayProof
			\qquad
				\bottomAlignProof
						\AxiomC{$\neg(A \land B)$}
							\AxiomC{$[\neg A]^{i}$} \noLine
						\UnaryInfC{$\vdots$} \noLine
						\UnaryInfC{$C$}
							\AxiomC{$[\neg B]^{i}$} \noLine
						\UnaryInfC{$\vdots$} \noLine
						\UnaryInfC{$C$}
					\RightLabel{$\neg \land E^{i}$}
					\TrinaryInfC{$C$}
				\DisplayProof

		\end{center}

		\vspace{0.25cm}

		\begin{center}

				\bottomAlignProof
						\AxiomC{$\neg A$}
						\AxiomC{$\neg B$}
					\RightLabel{$\neg \lor I$}
					\BinaryInfC{$\neg(A \lor B)$}
				\DisplayProof
			\qquad
				\bottomAlignProof
						\AxiomC{$\neg(A \lor B)$}
					\RightLabel{$\neg \lor E$}
					\UnaryInfC{$\neg A$}
				\DisplayProof
			$\!\!$
				\bottomAlignProof
						\AxiomC{$\neg(A \lor B)$}
					\UnaryInfC{$\neg B$}
				\DisplayProof
			\qquad
				\bottomAlignProof
						\AxiomC{$A$}
					\RightLabel{$DN$}
					\UnaryInfC{$\neg \neg A$}
				\DisplayProof
			$\!\!$
				\bottomAlignProof
						\AxiomC{$\neg \neg A$}
					\UnaryInfC{$A$}
				\DisplayProof

		\end{center}

		\vspace{0.25cm}

		For every $P \in \mathcal{P}^{\oid}$ and $1 \leq i \leq m$:

		\begin{center}

				\bottomAlignProof
						\AxiomC{$Pc_{1}\dots c_{m}$}
					\RightLabel{$\D I$}
					\UnaryInfC{$\D c_{i}$}
				\DisplayProof
			$\!\!$
				\bottomAlignProof
						\AxiomC{$\neg Pc_{1}\dots c_{m}$}
					\UnaryInfC{$\D c_{i}$}
				\DisplayProof
			\qquad
				\bottomAlignProof
						\AxiomC{$\D c$}
						\AxiomC{$\neg \D c$}
					\RightLabel{$PEX_{\D}$}
					\BinaryInfC{$B$}
				\DisplayProof
			\qquad
				\bottomAlignProof
						\AxiomC{}
					\RightLabel{$PEM_{\D}$}
					\UnaryInfC{$\D c \lor \neg \D c$}
				\DisplayProof

		\end{center}

		\vspace{0.25cm}

		\begin{center}

				\bottomAlignProof
						\AxiomC{$[\D c]^{i}$} \noLine
						\UnaryInfC{$\vdots$} \noLine
						\UnaryInfC{$A(c/x)$}
					\RightLabel{$\forall I^{i}_{\D}$}
					\UnaryInfC{$\forall xA$}
				\DisplayProof
			\qquad
				\bottomAlignProof
						\AxiomC{$\forall xA$}
						\AxiomC{$\D c$}
					\RightLabel{$\forall E_{\D}$}
					\BinaryInfC{$A(c/x)$}
				\DisplayProof
			\qquad
				\bottomAlignProof
						\AxiomC{$A(c/x)$}
						\AxiomC{$\D c$}
					\RightLabel{$\exists I_{\D}$}
					\BinaryInfC{$\exists xA$}
				\DisplayProof
			\qquad
				\bottomAlignProof
						\AxiomC{$\exists xA$}
							\AxiomC{$[A(c/x),\D c]^{i}$} \noLine
						\UnaryInfC{$\vdots$} \noLine
						\UnaryInfC{$B$}
					\RightLabel{$\exists E^{i}_{\D}$}
					\BinaryInfC{$B$}
				\DisplayProof

		\end{center}

		\vspace{0.25cm}

		\begin{center}

				\bottomAlignProof
						\AxiomC{$\neg A(c/x)$}
						\AxiomC{$\D c$}
					\RightLabel{$\neg \forall I_{\D}$}
					\BinaryInfC{$\neg \forall xA$}
				\DisplayProof
			\qquad
				\bottomAlignProof
						\AxiomC{$\neg \forall xA$}
							\AxiomC{$[\neg A(c/x), \D c]^{i}$} \noLine
						\UnaryInfC{$\vdots$} \noLine
						\UnaryInfC{$B$}
					\RightLabel{$\neg \forall E^{i}_{\D}$}
					\BinaryInfC{$B$}
				\DisplayProof

		\end{center}

		\vspace{0.25cm}

		\begin{center}

				\bottomAlignProof
						\AxiomC{$[\D c]^{i}$} \noLine
						\UnaryInfC{$\vdots$} \noLine
						\UnaryInfC{$\neg A(c/x)$}
					\RightLabel{$\neg \exists I^{i}_{\D}$}
					\UnaryInfC{$\neg \exists xA$}
				\DisplayProof
			\qquad
				\bottomAlignProof
						\AxiomC{$\neg \exists xA$}
						\AxiomC{$\D c$}
					\RightLabel{$\neg \exists E_{\D}$}
					\BinaryInfC{$\neg A(c/x)$}
				\DisplayProof

		\end{center}

		\vspace{0.25cm}

		\begin{center}

				\bottomAlignProof
						\AxiomC{$\D c$}
					\RightLabel{$\oid\!I$}
					\UnaryInfC{$c \oid c$}
				\DisplayProof
			\qquad
				\bottomAlignProof
						\AxiomC{$A(c_{1}/x)$}
						\AxiomC{$c_{1} \oid c_{2}$}
					\RightLabel{$\oid\!E$}
					\BinaryInfC{$A(c_{2}/x)$}
				\DisplayProof
			\qquad
				\bottomAlignProof
						\AxiomC{$A(c_{1}/x)$}
						\AxiomC{$\neg A(c_{2}/x)$}
					\RightLabel{$\noid\!I$}
					\BinaryInfC{$c_{1} \noid c_{2}$}
				\DisplayProof

		\end{center}

		\vspace{0.25cm}

		}

	\noi \textbf{Restrictions:} In $\forall I_{\D}$ ($\neg \exists I_{\D})$, $c$ must
	not occur in $A$, nor in any hypothesis on which $A(c/x)$ ($\neg
	A(c/x)$) depends, except $\D c$. In $\exists E_{\D}$ ($\neg \forall E_{\D}$), $c$
	must occur neither in $A$ or $B$, nor in any hypothesis on which $B$
	depends, except $A(c/x)$ ($\neg A(c/x)$) and $\D c$. In $\oid\!E$, $A$ is
	either an atomic formula or a negated atomic formula. In $\noid\!I$, $A$
	is an atomic formula in which $x$ is free.

	\end{df}

Given a first-order language \sig\ and $\Gamma \cup \{A\} \subseteq Sent(\sig)$,
the definition of a deduction of $A$ from $\Gamma$ in \ffde\ is the usual one
(see, e.g., \cite[Ch.~2]{Troelstra&vanDalen1988}). We shall use the notation
$\Gamma \vdash_{\sig} A$ to express that there is a deduction of $A$ from the
premises in $\Gamma$ in \ffde, omitting the subscript `$\sig$' when there is no
risk of confusion.

Note that rule $\oid\!I$ is equivalent to:

	\begin{center}

		\bottomAlignProof
				\AxiomC{}
			\UnaryInfC{$\forall x(x \oid x)$}
		\DisplayProof

	\end{center}

\noi as it is shown by the following deductions:

	\begin{center}

			\bottomAlignProof
						\AxiomC{$[\D c]^{1}$}
					\RightLabel{$\oid\!I$}
					\UnaryInfC{$c \oid c$}
				\RightLabel{$\forall I_{\D}^{1}$}
				\UnaryInfC{$\forall x(x \oid x)$}
			\DisplayProof
		\qquad
			\bottomAlignProof
					\AxiomC{}
					\UnaryInfC{$\forall x(x = x)$}
					\AxiomC{$\D c$}
				\RightLabel{$\forall E_{\D}$}
				\BinaryInfC{$c \oid c$}
			\DisplayProof

	\end{center}

\noi Note also that in $\noid\!I$ it is crucial that $x$ be free in $A$, for
otherwise $c_{1} \noid c_{2}$ would be derivable from \textit{any} pair of
contradictory atomic premises, for every $c_{1}, c_{2} \in \C$. Nonetheless, $c
\noid c$ \textit{is} derivable from any pair of such contradictions \textit{in
which} $c$ \textit{occurs}, e.g., $Pc$ and $\neg Pc$.

	\begin{rmrk} $\D c$ cannot be taken as an abbreviation of $\exists x(x
	\oid c)$ in \ffde, as it is customary in free logics. Given what was
	said about the intended interpretation of $\D$ and $\exists$, this
	should come as no surprise. However, this situation raises the question
	about what are the logical relationships between $\D c$ and $\exists x(x
	\oid c)$ in \ffde. The following deductions show that $\D c$ is
	equivalent to $\exists x(x \oid c)$:

		\begin{center}

				\bottomAlignProof
							\AxiomC{$\D c$}
						\RightLabel{$\oid\!I$}
						\UnaryInfC{$c \oid c$}
						\AxiomC{$\D c$}
					\RightLabel{$\exists I_{\D}$}
					\BinaryInfC{$\exists x(x \oid c)$}
				\DisplayProof
			\qquad
				\bottomAlignProof
					\AxiomC{$\exists x(x \oid c)$}
							\AxiomC{$[\D c']^{1}$}
							\AxiomC{$[c' \oid c]^{1}$}
						\RightLabel{$\oid\!E$}
						\BinaryInfC{$\D c$}
					\RightLabel{$\exists E_{\D}^{1}$}
					\BinaryInfC{$\D c$}
				\DisplayProof

		\end{center}

	\noi Morever, $\neg \exists x(x \oid c)$ is derivable from $\neg \D c$
	in \ffde:

		\begin{center}

			\bottomAlignProof
						\AxiomC{$[\D c']^{1}$}
						\AxiomC{$\neg \D c$}
					\RightLabel{$\noid\!I$}
					\BinaryInfC{$c' \noid c$}
				\RightLabel{$\neg \exists I_{\D}^{1}$}
				\UnaryInfC{$\neg \exists x(x \oid c)$}
			\DisplayProof

		\end{center}

	\noi Hence, $\exists x(x \oid c) \lor \neg \exists x(x \oid c)$ is a
	theorem of \ffde\ (by rule $PEM_{\D}$). However, $\neg \D c$ is not
	derivable from $\neg \exists x(x \oid c)$, since there are non-trivial
	models in which both $\exists x(x \oid c)$ and $\neg \exists x(x \oid
	c)$ hold (see Remark \ref{rmrk:ineqE} below). As a result, $\exists x(x
	\oid c)$ and $\neg \exists x(x \oid c)$ do not lead to triviality in
	\ffde, even though $\D c$ and $\neg \D c$ do.

	\end{rmrk}

The following proposition generalizes rules $\oid\!E$ and $\noid\!I$ to non-atomic
formulas:

	\begin{prop}\label{prop:igeneralization} Let $c_{1},c_{2} \in \C$ and
	$\Gamma \subseteq Sent(\sig)$ and suppose that $A \in Form(\sig)$ has at
	most $x$ free. Then:

		\begin{itemize}\setl

		\item[(1)] $\Gamma, A(c_{1}/x), c_{1} \oid c_{2} \vdash
		A(c_{2}/x)$;

		\item[(2)] If $x$ is free in every subformula of $A$, then
		$\Gamma, A(c_{1}/x), \neg A(c_{2}/x) \vdash c_{1} \noid c_{2}$.

		\end{itemize}

	\end{prop}

	\begin{pr} The proof of (1) straightforward, and so we shall only prove
	(2). Note first that $c_{2} \noid c_{1} \vdash c_{1} \noid c_{2}$:

		\begin{center}

			\bottomAlignProof
							\AxiomC{$c_{2} \noid c_{1}$}
						\RightLabel{$\D I$}
						\UnaryInfC{$\D c_{1}$}
					\RightLabel{$I \oid$}
					\UnaryInfC{$c_{1} \oid c_{1}$}
					\AxiomC{$c_{2} \noid c_{1}$}
				\RightLabel{$\noid\!I$}
				\BinaryInfC{$c_{1} \noid c_{2}$}
			\DisplayProof

		\end{center}

	\noi Let us call the associated rule $Sym\!\noid$. The proof of (2) runs
	by induction on the length of $A$. If $A$ is an atomic formula, then the
	result follows immediately from rule $\noid\!I$. If $A$ is $\neg B$ and
	$B$ is an atomic formula, then:

		\begin{center}

			\bottomAlignProof
							\AxiomC{$\neg \neg B(c_{2}/x)$}
						\RightLabel{$DN$}
						\UnaryInfC{$B(c_{2}/x)$}
						\AxiomC{$\neg B(c_{1}/x)$}
					\RightLabel{$\noid\!I$}
					\BinaryInfC{$c_{2} \noid c_{1}$}
				\RightLabel{$Sym\!\noid$}
				\UnaryInfC{$c_{1} \noid c_{2}$}
			\DisplayProof

		\end{center}

	\noi If $A$ is $\neg \neg B$, then:

		\begin{center}

			\bottomAlignProof
						\AxiomC{$\neg \neg B(c_{1}/x)$}
					\RightLabel{$DN$}
					\UnaryInfC{$B(c_{1}/x)$}
						\AxiomC{$\neg \neg \neg B(c_{2}/x)$}
					\RightLabel{$DN$}
					\UnaryInfC{$\neg B(c_{2}/x)$}
				\RightLabel{$IH$}
				\BinaryInfC{$c_{1} \noid c_{2}$}
			\DisplayProof

		\end{center}

	\noi If $A$ is $\neg (B \land C)$, then:

		\begin{center}

			\bottomAlignProof
							\AxiomC{$\neg \neg(B(c_{2}/x) \land C(c_{2}/x))$}
						\RightLabel{$DN$}
						\UnaryInfC{$B(c_{2}/x) \land C(c_{2}/x)$}
						\AxiomC{$\neg(B_{1}(c_{1}/x) \land C(c_{1}/x))$}
					\RightLabel{$IH$}
					\BinaryInfC{$c_{2} \noid c_{1}$}
				\RightLabel{$Sym\!\noid$}
				\UnaryInfC{$c_{1} \noid c_{2}$}
			\DisplayProof

		\end{center}

	\noi Let $A$ be $\forall yB$. Note that $(\forall yB)(c_{1}/x) = \forall
	yB(c_{1}/x)$ and $(\forall yB)(c_{2}/x) = \forall yB(c_{2}/x)$, since
	$x$ is free in $A$, and so $x \neq y$. Note also that $B(c_{1}/x)(c/y) =
	B(c/y)(c_{1}/x)$ and $B(c_{2}/x)(c/y) = B(c/y)(c_{2}/x)$, for every $c
	\in \mathcal{C}$. Thus:

		\begin{center}

			\bottomAlignProof
					\AxiomC{$\neg \forall yB(c_{2}/x)$}
							\AxiomC{$\forall yB(c_{1}/x)$}
							\AxiomC{$[\D c]^{1}$}
						\RightLabel{$\forall E_{\D}$}
						\BinaryInfC{$B(c_{1}/x)(c/y)$}
						\UnaryInfC{$B(c/y)(c_{1}/x)$}
						\AxiomC{$[\neg B(c_{2}/x)(c/y)]^{1}$}
						\UnaryInfC{$[\neg B(c/y)(c_{2}/x)]^{1}$}
					\RightLabel{$IH$}
					\BinaryInfC{$c_{1} \noid c_{2}$}
				\RightLabel{$\neg \forall E^{1}_{\D}$}
				\BinaryInfC{$c_{1} \noid c_{2}$}
			\DisplayProof

		\end{center}

	\noi The remaining cases are similar to the ones above.

	\end{pr}

Note that we would be unable prove Proposition \ref{prop:igeneralization}(2) if
the requirement that $x$ be free in every subformula of $A$ were missing, for
otherwise we would be prevented from applying the indutive hypothesis to the
immediate subformulas of, say, $B \land C$ (see Remark
\ref{rmrk:igeneralization} for a counter-example).

\subsection{Semantics}\label{sec:semantics}

A model for \ffde\ is composed of a set $W$ of \textit{stages}, a binary
relation $\leq$ among the elements of $W$, an interpretation function $I$, a
function $d$ that assigns sets to the elements of $W$, and a pair of functions
$i_{+}$ and $i_{-}$ on $W$. For each stage $w \in W$, $d(w)$ is the
\textit{domain} of $w$, that is, the range of the quantifiers at $w$, while
$i_{+}(w)$ and $i_{-}(w)$ determine the extension and anti-extension of the
identity predicate $\oid$ at $w$.

	\begin{df}\label{df:kripkemodel} Let $\sig = \langle \C, \mathcal{P}
	\rangle$ be a first-order language. A (\textit{Kripke}
	\ffde-)\textit{model} \K\ for \sig\ is a structure $\langle W, \leq, d,
	I, i_{+}, i_{-} \rangle$ such that: (i) $W$ is a non-empty set of
	\textit{stages}; (ii) $\leq$ (\textit{accessibility relation}) is a
	pre-order on $W$; (iii) $d$ is a function on $W$ such that $d(w)$ is a
	set (the \textit{domain} of $w$) satisfying the following condition: if
	$w \leq w'$, then $d(w) \subseteq d(w')$; (iv) $I$
	(\textit{interpretation function}) is a partial function on $(\C \cup
	\mathcal{P}) \times W$ such that:

		\begin{itemize}\setl

		\item For every individual constant $c \in \C$, if $I(c,w)$ is
		defined, then $I(c, w) \in d(w)$;

		\item For every individual constant $c \in \C$, if $I(c,w)$ is
		defined, then $I(c, w) = I(c, w')$, for every $w' \geq w$;

		\item For every $m$-ary predicate letter $P \in \mathcal{P}$,
		$I(P, w)$ is the pair $\langle P^{w}_{+}, P^{w}_{-} \rangle$
		such that $P^{w}_{+}, P^{w}_{-} \subseteq d(w)^{m}$ and if $w
		\leq w'$, then $P^{w}_{+} \subseteq P^{w'}_{+}$ and $P^{w}_{-}
		\subseteq P^{w'}_{-}$.

		\end{itemize}

	\noi (v) $i_{+}$ and $i_{-}$ are functions on $W$ such that $i_{+}(w)$
	and $i_{-}(w)$ are binary relations on $d(w)$ that satisfy the following
	conditions:

		\begin{itemize}

		\item[(a)] $i_{+}(w)$ is a congruence relation on $d(w)$. That is,
		$i_{+}(w)$ is an equivalence relation such that for every $P \in
		\mathcal{P}$ and $a_{1},\dots,a_{m},a_{1}',\dots,a_{m}' \in
		d(w)$, if $\langle a_{j}, a_{j}' \rangle \in i_{+}(w)$, for
		every $1 \leq j \leq m$, then $\langle a_{1},\dots,a_{m} \rangle
		\in P^{w}_{+}$ iff $\langle a_{1}',\dots,a_{m}' \rangle \in
		P^{w}_{+}$, and $\langle a_{1},\dots,a_{m} \rangle \in P^{w}_{-}$
		iff $\langle a_{1}',\dots,a_{m}' \rangle \in P^{w}_{-}$;

		\item[(b)] For every $a, a' \in d(w)$, if there is an
		$m$-ary predicate $P \in \mathcal{P}^{\oid}$ and an $m$-tuple
		$\langle \dots a \dots \rangle$ of elements of $d(w)$ such that
		$\langle \dots a \dots \rangle \in P^{w}_{+}$ and $\langle \dots
		a' \dots \rangle \in P^{w}_{-}$, then $\langle a, a' \rangle
		\in i_{-}(w)$ and $\langle a', a \rangle \in i_{-}(w)$.

		\item[(c)] For every $w, w' \in W$, if $w \leq w'$, then
		$i_{+}(w) \subseteq i_{+}(w')$ and $i_{-}(w) \subseteq
		i_{-}(w')$.

		\end{itemize}

	\end{df}

Note that the domain of a stage is not required to be non-empty, but must
include the domains of all $\leq$-preceding stages. Note also that according to
clause (iv), $I$ assigns the same individual to $c \in \C$ across all stages.
Thus, there can be no $a$ such that $I(c, w) = a$ but $a \notin d(w')$, for some
$w' \geq w$, and the domains of two stages can only differ with respect to
``unamed'' individuals. The interpretation of the predicate letters, on the
other hand, may vary across stages: for each $w \in W$, $P \in \mathcal{P}$ is
assigned both an \textit{extension} $P^{w}_{+}$ and an \textit{anti-extension}
$P^{w}_{-}$ such that $P^{w}_{+}$ represents the set of individuals \textit{at}
$w$ that are ``known'' to satisfy $P$, while $P^{w}_{-}$ represents the set of
individuals \textit{at} $w$ that are ``known'' to satisfy not-$P$. Although
$P^{w}_{+} \cup P^{w}_{-}$ must be a subset of $d(w)$, it may not be identical
to $d(w)$. Moreover, the intersection of $P^{w}_{+}$ and $P^{w}_{-}$ may be
non-empty. According to the intended intuitive interpretation of \ffde, whenever
$P^{w}_{+} \cup P^{w}_{-} \neq d(w)$ ($P^{w}_{+} \cap P^{w}_{-} \neq
\emptyset$), then there will be incomplete (contradictory) information involving
$P$. The extension $i_{+}(w)$ and the anti-extension $i_{-}(w)$ of $\oid$ at $w$
behave as $P^{w}_{+}$ and $P^{w}_{-}$, except that they are subject to the
restrictions discussed towards the end of the last section, viz., conditions (a)
and (b). In particular, (b) states that if the extension and the anti-extension
of an $m$-ary predicate letter $P$ at $w$ include respectively two $m$-tuples
$\langle \dots a \dots \rangle$ and $\langle \dots a' \dots \rangle$ that differ
at most with respect to the terms $a$ and $a'$ (which are assumed to occupy the
same position in each tuple), then the pairs $\langle a, a' \rangle$ and
$\langle a, a' \rangle$ must belong to $i_{-}(w)$. This means that the model
``sees'' these two individuals as distinct at $w$. The requirement that both
pairs $\langle a, a' \rangle$ and $\langle a', a \rangle$ belong to $i_{-}(w)$
is meant to ensure that $\noid$ is indeed a symmetric relation, as one might
expect it to be.

Since we will adopt a substitutional interpretation of the quantifiers, the
following definition will come in handy:

	\begin{df}\label{df:diagram} Let $\sig = \langle \C, \mathcal{P}
	\rangle$ be a first-order language and let \K\ be a Kripke model for
	\sig. The \textit{diagram language} $\sig_{\K}$ of \K\ is the pair
	$\langle \C_{\K}, \mathcal{P} \rangle$ such that $\C_{\K} = \C \cup
	\{\ov{a}: a \in \bigcup_{w\in W} d(w)\}$. We will use the notation
	$\wh{\K}$ to refer to the Kripke model $\langle W, \leq, d, \wh{I},
	i_{+}, i_{-} \rangle$ for $\sig_{\K}$ that is just like \K\ except that
	$\wh{I}$ is such that $\wh{I}(\ov{a}, w) = a$ if $a \in d(w)$, and
	undefined otherwise.

	\end{df}

\noi Thus, the diagram language $\sig_{\K}$ of \sig\ results by adding a new
individual constant $\ov{a}$ for each $a$ in the domain of at least one stage of
\K, while the corresponding model $\wh{\K}$ is such that its interpretation
$\wh{I}$ function assigns $a$ to $\ov{a}$ at each stage $a$ at which ``exists''.

	\begin{df}\label{df:valuation} Let \sig\ be a first-order language and
	let \K\ be a Kripke model for \sig. The \textit{valuation function
	induced by} \K\ is the mapping $v\!: Sent(\sig_{\K}) \times W
	\longrightarrow \{0,1\}$ satisfying the following conditions:

		\begin{itemize}\setl

		\item[($v$1)] $v(\D c,w) = 1$ iff $\wh{I}(c,w)$ is defined.

		\item[($v$2)] $v(\neg \D c,w) = 1$ iff $\wh{I}(c,w)$ is undefined.

		\item[($v$3)] $v(c_{1} \oid c_{2}, w) = 1$ iff $\wh{I}(c_{1},w)$
		and $\wh{I}(c_{2},w)$ are defined and $\langle \wh{I}(c_{1},w),
		\wh{I}(c_{1},w) \rangle \in i_{+}(w)$;

		\item[($v$4)] $v(c_{1} \noid c_{2}, w) = 1$ iff
		$\wh{I}(c_{1},w)$ and $\wh{I}(c_{2},w)$ are defined and $\langle
		\wh{I}(c_{1},w), \wh{I}(c_{1},w) \rangle \in i_{-}(w)$;

		\item[($v$5)] For every $m$-ary predicate letter $P \in
		\mathcal{P}$, $v(Pc_{1}\dots c_{m},w) = 1$ iff $\wh{I}(c_i,w)$
		is defined, for every $1 \leq i \leq m$, and $\langle
		\wh{I}(c_{1},w),\dots,\wh{I}(c_{m},w) \rangle \in P^{w}_{+}$;

		\item[($v$6)] For every $m$-ary predicate letter $P \in
		\mathcal{P}$, $v(\neg Pc_{1}\dots c_{m},w) = 1$ if $\wh{I}(c_i,
		w)$ is defined, for every $1 \leq i \leq m$, and $\langle
		\wh{I}(c_{1},w),\dots,\wh{I}(c_{m},w) \rangle \in P^{w}_{-}$;

		\item[($v$7)] $v(A \land B,w) = 1$ iff $v(A,w) = 1$ and $v(B,w)
		= 1$;

		\item[($v$8)] $v(A \lor B,w) = 1$ iff $v(A,w) = 1$ or $v(B,w) =
		1$;

		\item[($v$9)] $v(\neg (A \land B),w) = 1$ iff $v(\neg A,w) = 1$
		or $v(\neg B,w) = 1$;

		\item[($v$10)] $v(\neg (A \lor B),w) = 1$ iff $v(\neg A,w) = 1$
		and $v(\neg B,w) = 1$;

		\item[($v$11)] $v(\neg \neg A,w) = 1$ iff $v(A,w) = 1$;

		\item[($v$12)] $v(\forall xA,w) = 1$ iff for every $w' \geq w$,
		$v(A(\ov{a}/x),w') = 1$, for every $a \in d(w')$;

		\item[($v$13)] $v(\exists xA,w) = 1$ iff $v(A(\ov{a}/x), w) = 1$,
		for some $a \in d(w)$;

		\item[($v$14)] $v(\neg \forall xA,w) = 1$ iff $v(\neg
		A(\ov{a}/x),w) = 1$, for some $a \in d(w)$;

		\item[($v$15)] $v(\neg \exists xA,w) = 1$ iff for every $w' \geq
		w$, $v(\neg A(\ov{a}/x), w') = 1$, for every $a \in d(w')$.

		\end{itemize}

	\end{df}

\noi A sentence $A \in Sent(\sig)$ \textit{holds} at a stage $w$ of \K\ ($\K, w
\vDash A$) if $v(A,w) = 1$, and a set of sentences $\Gamma \subseteq Sent(\sig)$
\textit{holds} at $w$ ($\K,w \vDash \Gamma$) if $\K, w \vDash B$, for every $B
\in \Gamma$. Given $\Gamma \cup \{A\} \subseteq Sent(\sig)$, $A$ is a
\textit{semantic consequence} of $\Gamma$ ($\Gamma \vDash A$) if and only if for
every Kripke model \K\ and every stage $w$ of \K, $\K,w \vDash \Gamma$ only if
$\K,w \vDash A$.

We have already seen that from a proof-theoretical viewpoint \ffde\ is indeed a
free logic, for it replaces the usual $\forall$-elimination and
$\exists$-introduction rules (as well as the usual $\neg \forall$-introduction
and $\neg \exists$-elimination rules) by weaker versions thereof that require
the relevant constant $c$ to satisfy the \ffde\ analogue of the existence
predicate, viz., \D. The new versions of these rules, with \E\ replacing \D, are
often present in standard formulations of natural deduction systems for several
free logics. Also, the replacement of the usual $\forall$-introduction and
$\exists$-elimination rules (as well as the $\neg \forall$-elimination and $\neg
\exists$-introduction rules) by respectively $\forall I_{\D}$ and $\exists
E_{\D}$ ($\neg \forall E_{\D}$ and $\neg \exists I_{\D}$) turns it into an
inclusive logic, and so \ffde\ is universally free. From a semantic perspective,
the inclusiveness of \ffde\ is of course due to the fact that $d(w)$ and $d(w')$
in clauses ($v$12)-($v$15) of Definition \ref{df:valuation} are not required to
denote non-empty sets. Note that $\forall xA$ and $\neg \exists xA$ vacuously
hold at every empty stage $w$ of $\K$, provided that $d(w') = \emptyset$, for
every $w' \geq w$. Also, $\exists xA$ and $\neg \forall xA$ do not hold at any
empty stage, irrespective of whether $d(w') \neq \emptyset$, for some $w' \geq
w$. Moreover, since a constant $c$ may be empty at $w$, there is no guarantee
that $A(c/x)$ holds at $w$ given that $A$ holds of every element of $d(w')$, for
every $w' \geq w$. In particular, if $I(c,w)$ is undefined and $A$ is atomic,
then $v(A(c/x), w) = 0$, even though $v(A(\ov{a}/x), w') = 1$, for every $w'
\geq w$ and every $a \in d(w')$. However, if $\D c$ holds at $w$, then $I(c,w)$
is defined, which means that $I(c,w) \in d(w')$, for every $w' \geq w$ (by
Definition \ref{df:kripkemodel}(iv)). Hence, if $\forall xA$ holds at $w$, we
may expect $A(c/x)$ to do so as well.

	\begin{rmrk}

	According to Definition \ref{df:valuation}($v$3)-($v$6), every atomic
	sentence receives the value $0$ at a stage $w$ whenever it contains an
	occurrence of an empty constant (and so does the negation of any such
	sentence not of the form $\neg \D c$). This fact raises the question of
	whether \ffde\ should be classified as a \textit{negative} or as a
	\textit{neutral free logic}. In a negative free logic every atomic
	sentence in which an empty term occurs is invariably false, while in a
	neutral logic every such sentence not of the form $\E t$ is neither true
	nor false (by contrast, in a \textit{positive free logic} some atomic
	sentences with empty terms may be true). Since the standard definitions
	of these labels are formulated in terms alethic values, they cannot be
	readly applied to \ffde, and the answer to the present question depends
	on how one chooses to extend the above classification to logics with
	non-alethic semantic values. Specifically, if one takes a negative
	logic to be any free logic that assigns a non-designated value to every
	atomic formula in which an empty term occurs, then \ffde\ is definitely
	a negative logic. If, on the other hand, one takes into account that
	according to the informal interpretation of the values $1$ and $0$,
	$v(A,w) = v(\neg A, w) = 0$ means that there is no information about $A$
	at $w$, it appears to be more appropriate to regard \ffde\ as neutral
	logic. Since both ways of classifying \ffde\ seem equally correct and
	nothing that will be said below will depend on these labels, we won't
	take sides concerning this issue.

	\end{rmrk}

The following proposition states that if $A$ holds at a stage $w$ and $w' \geq
w$, then $A$ also holds at $w'$. Hence, the value $1$ assigned by the valuation
induced by a model persists across $\leq$-related stages. According to the
intuitive interpretation of \ffde\ as an information-based logic, this means
that once a certain piece of information has been added to a database it can no
longer be removed.

	\begin{prop}\label{prop:persistence} Let $A \in Sent(\sig)$ and let \K\
	be a Kripke model for \sig. For every stages $w$ and $w'$ of \K, if $\K,
	w \vDash A$ and $w \leq w'$, then $\K, w' \vDash A$.

	\end{prop}

	\begin{pr} The result follows by a straightforward induction on the
	complexity of $A$.

	\end{pr}

It is worth observing at this point that Proposition \ref{prop:persistence} is
in line with Belnap's idea of how the information states change over time:

	\begin{quote}

	When an atomic formula is entered into the computer as either affirmed
	or denied, the computer modifies its current set-up by adding a `told
	True' or `told False' according as the formula was affirmed or denied;
	it does not subtract any information it already has for that is the
	whole point of what we are up to. In other words, if $p$ is affirmed, it
	marks $p$ with \textbf{T} if $p$ were previously marked with
	\textbf{None}; with \textbf{Both} if $p$ were previously marked with
	\textbf{False} [sic.] \cite[p.~36]{Belnap1977}.

	\end{quote}

\noi Indeed, according to the interpretation of the values $1$ and $0$ adopted
here, the persistence property disallows ``information to be subtracted'' from
the computer, as long as we understand the four scenarios described by Belnap in
the following way (where $v_{4}$ is a four-valued valuation):

	\begin{itemize}\setl

	\item[] $v_4(A) = \textbf{T}$ iff $v(A)=1 $ and $ v(\neg A)=0$,
	\item[] $v_4(A) = \textbf{F}$ iff $v(A)=0$ and $ v(\neg A)=1$,
	\item[] $v_4(A) = \textbf{None}$ iff $v(A)=0 $ and $ v(\neg A)=0$,
	\item[] $v_4(A) = \textbf{Both}$ iff $v(A)=1$ and $v(\neg A)=1$.

	\end{itemize}

\noi As expected, these two ways of formulating the semantics of \fde\ are
equivalent (see \cite[Sect.~2.2]{Antunesetal2022}), and these are equivalent to
the relational semantics proposed by Dunn \citep{Dunn1976}.

In the passage above Belnap also observes that the value \textbf{N} can change
to \textbf{T}, while \textbf{F} can change to \textbf{B}. In view of the
translation of Belnap's semantics into the bivalued semantics adopted here, this
corresponds to allowing the value $0$ to change to $1$, which, in turn,
represents the fact that positive (negative) information has been added to a
database: when the value of $A$ ($\neg A$) changes from $0$ to $1$, this
corresponds, in the four-valued semantics, to the value \textbf{N} changing to
\textbf{T}, provided that the value of $\neg A$ ($A$) is also $0$; and to the
value \textbf{F} changing to \textbf{B}, provided that the value of $\neg A$
($A$) is $1$. As a way of an example, consider the model $\K = \langle W, \leq,
d, I, i_{+}, i_{-} \rangle$ such that:

	\begin{itemize}\setl

	\item $W = \{w_{1}, w_{2}\}$; $\leq \ = \{\langle w_{1}, w_{1}
	\rangle, \langle w_{2}, w_{2} \rangle, \langle w_{1}, w_{2}
	\rangle\}$;

	\item $d(w_{1}) = d(w_{2}) = \{a\}$;

	\item $I(c) = a$, $P^{w_{1}}_{+} = P^{w_{1}}_{-} = \emptyset$,
	$P^{w_{2}}_{+} = \{a\}$, and $P^{w_{2}}_{-} = \emptyset$;

	\item $Q^{w_{1}}_{-} = \emptyset$, $Q^{w_{1}}_{+} = Q^{w_{2}}_{+} =
	\{a\} = Q^{w_{2}}_{-} = \{a\}$;

	\item $i_{+}(w_{1}) = i_{+}(w_{2}) = \{\langle a, a \rangle \}$,
	$i_{-}(w_{1}) = \emptyset$, and $i_{-}(w_{2}) = \{\langle a, a
	\rangle\}$.

	\end{itemize}

\noi In this model $v(Pc,w_{1}) = v(\neg Pc,w_{1}) = 0$, while $v(Pc,w_{2}) =
1$, which corresponds to the value of $Pc$ changing from \textbf{N} to
\textbf{T} in Belnap's semantics. As for $Qc$, notice that $v(Qc,w_{1}) = 1$,
while $v(\neg Qc,w_{1}) = 0$ and $v(Qc,w_{2}) = v(\neg Qc,w_{2}) = 1$, which
corresponds to the value of $\neg Qc$ changing from \textbf{T} to \textbf{B}.
According to the intended interpretation of \ffde, positive information $Pc$ and
negative information $Qc$ have been added to the database as it evolves from
$w_{1}$ to $w_{2}$. In the latter case, though, notice that positive information
$Qc$ was already available at $w_{1}$, and so the database becomes contradictory
at $w_{2}$.

Let us now show that the (inferential version of the) constant domain axiom does
not hold in \ffde.

	\begin{prop}\label{prop:constantdomain} There is a model \K\ and a stage
	$w$ of \K\ such that $\K, w \vDash \forall x(B \lor A)$ and $\K, w
	\nvDash B \lor \forall xA$, for some $A$ and $B$ such that $A$ has at
	most $x$ free and $B$ is a sentence.

	\end{prop}

	\begin{pr} Let $\K = \langle W, \leq, d, I, i_{+}, i_{-} \rangle$ be
	such that:

		\begin{itemize}\setl

		\item $W = \{w_{1}, w_{2}\}$; $\leq \ = \{\langle w_{1}, w_{1}
		\rangle, \langle w_{2}, w_{2} \rangle, \langle w_{1}, w_{2}
		\rangle\}$;

		\item $d(w_{1}) = \{a_{1}\}$ and $d(w_{2}) = \{a_{1},a_{2}\}$;
		$P^{w_{1}}_{+} = P^{w_{2}}_{+} = \{a_{1}\}$, $P^{w_{1}}_{-} =
		\emptyset$, and $P^{w_{2}}_{-} = \{a_{2}\}$;

		\item $i_{+}(w_{1}) = \{\langle a_{1}, a_{1} \rangle \}$,
		$i(w_{2}) = \{\langle a_{1},a_{1} \rangle, \langle a_{2}, a_{2}
		\rangle \}$, $i_{-}(w_{1}) = \emptyset$ and $i_{-}(w_{2}) =
		\{\langle a_{1}, a_{2} \rangle, \langle a_{2}, a_{1} \rangle \}$.

		\end{itemize}

	\noi Let $B = \exists y \neg Py$. Clearly, $\K, w_{1} \nvDash B$ and
	$\K, w_{2} \vDash B$. Now, since $\K, w_{1} \vDash B \lor P\ov{a_{1}}$,
	$\K, w_{2} \vDash B \lor P\ov{a_{1}}$, and $\K, w_{2} \vDash B \lor
	P\ov{a_{2}}$, it follows that $\K, w_{1} \vDash \forall x(B \lor Px)$.
	But since $\K, w_{1} \nvDash B$ and $\K, w_{1} \nvDash \forall xPx$,
	$\K, w_{1} \nvDash B \lor \forall xPx$.

	\end{pr}

In Remark \ref{rmrk:ineqE} below we show that even though $\D c$ is equivalent
to $\exists x(x \oid c)$ and $\neg \D c$ entails $\neg \exists x(x \oid c)$ in
\ffde, $\neg Dc$ is not a logical consequence of $\neg \exists x(x \oid c)$. In
Remark \ref{rmrk:igeneralization} we present a counter-example to Proposition
\ref{prop:igeneralization}(2) when the condition that $x$ be free in every
subformula of $A$ is dropped.

	\begin{rmrk}\label{rmrk:ineqE} Let $\sig = \langle \C, \mathcal{P}
	\rangle$ be such that $\C = \{c\}$ and $\mathcal{P} = \emptyset$. Both
	$\exists x(x \oid c)$ and $\neg \exists x(x \oid c)$ hold in every stage
	of the model $\K = \langle W, \leq, d, I, i_{+}, i_{-} \rangle$ such
	that $W = \{w\}$, $\leq \ = \{\langle w, w \rangle \}$, $d(w) = \{1\}$,
	$I(c,w) = 1$, $i_{+}(w) = \{\langle 1, 1 \rangle\}$, and $i_{-}(w) =
	\{\langle 1, 1 \rangle\}$. Since $I(c)$ is defined, $\K, w \nvDash \neg
	\D c$. Thus, $\neg \exists x(x \oid c) \nvDash \neg \D c$.

	\end{rmrk}

	\begin{rmrk}\label{rmrk:igeneralization} Let $\sig = \langle \C,
	\mathcal{P} \rangle$ be such that $\C = \{c_{1}, c_{2}\}$ and
	$\mathcal{P} = \{P,Q\}$, where $P$ and $Q$ are unary. Let $\K = \langle
	W, \leq, d, I, i_{+}, i_{-} \rangle$ be an \sig-model such that $W =
	\{w\}$, $d(w) = \{1,2\}$, $i_{+}(w) = \{\langle 1, 1 \rangle, \langle 2,
	2 \rangle\}$, and $i_{-}(w) = \{\langle 2, 2 \rangle\}$. Suppose that
	$I(c_{1}, w) = 1$ and $I(c_{2},w) = 2$. Suppose further that $Q^{w}_{+}
	= \{2\}$ and $Q^{w}_{-} = \{2\}$. Thus, $\K, w \vDash Qc_{2}$ and $\K, w
	\vDash \neg Qc_{2}$. Letting $P^{w}_{+} = \{1\}$ and $P^{+}_{-} =
	\emptyset$, it then follows that $\K, w \vDash Pc_{1} \land Qc_{2}$ and
	$\K, w \nvDash \neg (Pc_{1} \land Qc_{2})$. Therefore, even though both
	sentences hold in $\K$ and $w$, $c_{1} \noid c_{1}$ does not.

	\end{rmrk}

\subsection{Soundness and Completeness}\label{sec:soundness_and_completeness}

We turn now to the task of proving soundness and completeness results for \ffde.
The proofs proceed along the lines of those in \cite{Rodrigues&Antunes2022}, for the logic
\qvletf, and in \cite{Antunes&Rodrigues2024}, for the variable domain (non-free) version
\qvfde\ of \ffde.

The following two technical lemmas are necessary for proving the soundness of
the quantifier rules of \ffde:

	\begin{lm}\label{lm:substitution_01} Let $\sig = \langle \mathcal{P}, \C
	\rangle$ be a first-order language such that $c \in \C$. Suppose that
	$A \in Form(\sig)$ has at most $x$ free. Let $\K = \langle W, \leq, d,
	I, i_{+}, i_{-} \rangle$ be a model and $w \in W$. Assume that $I(c, w)$
	is defined and let $I(c,w) = a \in d(w)$. Then $\K, w \vDash A(c/x)$ iff
	$\wh{\K}, w \vDash A(\ov{a}/x)$.

	\end{lm}

	\begin{lm}\label{lm:substitution_02} Let $\sig = \langle \mathcal{P},
	\C \rangle$ be a first-order language such that $c_{1},c_{2} \in \C$.
	Let $\K = \langle W, \leq, d, I, i_{+}, i_{-} \rangle$ be a model and $w
	\in W$. Assume that $I(c_{2},w)$ is defined and let $I(c_{2},w) = a \in
	d(w)$. Let $\K_{c_{2}/c_{1}} = \langle W, \leq, d, I_{c_{2}/c_{1}},
	i_{+}, i_{-} \rangle$ be such that $I_{c_{2}/c_{1}}$ agrees with $I$ in
	all elements of $\mathcal{P} \cup \mathcal{C}\setminus\{c_{1}\}$. Assume that
	$I_{c_{2}/c_{1}}(c_{1}, w') = a$, for every $w' \in W$ such that
	$I(c_{2},w')$ is defined ($I_{c_{2}/c_{1}}(c_{1},w')$ is undefined
	otherwise). If $A \in Form(\sig)$ has at most $x$ free and $c_{1}$ does
	not occur in $A$, then $\K_{c_{2}/c_{1}}, w \vDash A(c_{1}/x)$ if and
	only if $\K, w \vDash A(c_{2}/x)$.

	\end{lm}

The proofs of lemmas \ref{lm:substitution_01} and \ref{lm:substitution_02}
follow by straightforward inductions on the complexity of $A$.

	\begin{tr}\label{tr:soundness} \textbf{(Soundness)} Let $\Gamma \cup
	\{A\} \subseteq Sent(\sig)$. If $\Gamma \vdash_{\sig} A$, then $\Gamma
	\vDash A$.

	\end{tr}

	\begin{pr} The result is proven by induction on the number of nodes of a
	deduction of $A$ from $\Gamma$. For the sake to brevity, we shall only
	prove the soundness of a few rules.

		\begin{itemize}\setl

		\item $\D I$: Suppose that $\Gamma \vDash Pc_{1}\dots c_{m}$ and
		assume that $\K, w \vDash \Gamma$.\footnote{The proof of the
		soundness of the second version of $\D I$ is nearly identical.}
		Thus, $\K, w \vDash Pc_{1}\dots c_{n}$, and so $\wh{I}(c_{i},w)$
		is defined, for every $1 \leq i \leq m$. Hence, $\K, w \vDash \D
		c_{i}$, for every $1 \leq i \leq m$.

		\item $PEX_{\D}$: Suppose that $\Gamma \vDash \D c$ and $\Gamma
		\vdash \neg \D c$, and assume that $\K, w \vDash \Gamma$. Thus,
		$\K, w \vDash \D c$ and $\K, w \vDash \neg \D c$. By clauses
		($v1$) and ($v2$) of Definition \ref{df:valuation}, it follows
		that $\wh{I}(c,w)$ is both defined and undefined. Therefore,
		$\K, w \vDash B$.

		\item $\forall E_{\D}$: Suppose that $\Gamma \vDash \forall xA$
		and $\Gamma \vDash \D c$, and assume that $\K, w \vDash \Gamma$.
		Since $\K, w \vDash \D c$, $I(c,w)$ is defined. Let $I(c,w) = a
		\in d(w)$. Since $\K, w \vDash \forall xA$, $\wh{\K}, w \vDash
		A(\ov{a}/x)$. By Lemma \ref{lm:substitution_01}, it follows that
		$\K, w \vDash A(c/x)$.

		\item $\forall I_{\D}$: Suppose that $\Gamma, \D c \vDash
		A(c/x)$ and assume that $\K, w \vDash \Gamma$. Let $w' \geq w$.
		If $d(w') = \emptyset$, it vacuously follows that $\wh{\K}, w'
		\vDash A(\ov{a}/x)$, for every $a \in d(w')$. If $d(w') \neq
		\emptyset$, let $a$ be an arbitrary element of $d(w')$. By Lemma
		\ref{lm:substitution_02}, $\wh{\K}_{\ov{a}/c}, w' \vDash
		\Gamma$, since $\wh{\K}, w' \vDash \Gamma$ and $c$ does not
		occur in $\Gamma$. Moreover, $\wh{\K}_{\ov{a}/c}, w' \vDash \D
		c$, since $\wh{I}_{\ov{a}/c}(c,w') = a$. Hence,
		$\wh{\K}_{\ov{a}/c}, w' \vDash \Gamma \cup \{\D c\}$, and so
		$\wh{\K}_{\ov{a}/c}, w' \vDash A(c/x)$. By Lemma
		\ref{lm:substitution_02}, $\wh{\K}, w' \vDash A(\ov{a}/x)$. But
		because $a$ was arbitrary, it follows that $\K, w' \vDash
		A(\ov{a}/x)$, for every $a \in d(w')$. Since $w' \geq w$ was
		arbitrary, it then follows that $\K, w \vDash \forall xA$.

		\item $\exists E_{\D}$: Suppose that $\Gamma \vDash \exists xA$
		and $\Gamma, A(c/x), \D c \vDash B$, and assume that $\K, w
		\vDash \Gamma$. Thus, $\K, w \vDash \exists xA$, and so there is
		$a \in d(w)$ such that $\wh{\K}, w \vDash \{A(\ov{a}/x), \D
		\ov{a}\}$. By Lemma \ref{lm:substitution_02}, it follows that
		$\wh{\K}_{\ov{a}/c}, w \vDash \{A(c/x), \D c\}$ and
		$\wh{\K}_{\ov{a}/c}, w \vDash \Gamma$ (since $c$ does occurs
		neither in $A$ nor in $\Gamma$). Therefore, $\wh{\K}_{\ov{a}/c},
		w \vDash B$. But because $c$ does not occur in $B$, if follows
		by Lemma \ref{lm:substitution_02} that $\wh{\K}, w \vDash B$,
		and so $\K, w \vDash B$.

		\item $\noid\!I$: Let $A$ be an atomic formula in which $x$ is
		free and in which the $m$-ary predicate letter $P$ occurs.
		Suppose that $\Gamma \vDash A(c_{1}/x)$ and $\Gamma \vDash \neg
		A(c_{2}/x)$. Assume that $\K, w \vDash \Gamma$. Hence, $\K, w
		\vDash A(c_{1}/x)$ and $\K, w \vDash \neg A(c_{2}/x)$. By
		clauses ($v3$)-($v6$) of Definition \ref{df:valuation}, there
		are $m$-ary tuples $\langle \dots \wh{I}(c_{1},w) \dots \rangle$
		and $\langle \dots \wh{I}(c_{2},w) \dots \rangle$ such that
		$\langle \dots \wh{I}(c_{1},w) \dots \rangle \in P^{w}_{+}$ and
		$\langle \dots \wh{I}(c_{2},w) \dots \rangle \in P^{w}_{-}$. It
		follows from Definition \ref{df:kripkemodel}(v) that $\langle
		\wh{I}(c_{1},w), \wh{I}(c_{2},w) \rangle \in i_{-}(w)$, and so
		$\K, w \vDash c_{1} \noid c_{2}$.

		\end{itemize}

	\end{pr}

Let us now prove that \ffde\ is strongly complete with respect to the semantics
presented above. The proof will be divided into two main steps. In the first
step we show how any set of sentences $\Gamma$ that does not derive a certain
sentence $A$ can be extended to a set $\Delta$ that enjoys several desirable
properties while still not deriving $A$. The second step involves constructing a
canonical model \K\ all of whose stages are regular Henkin sets, and then
proving that although $\Gamma$ holds in \K, $\K \nvDash A$.

	\begin{df}\label{df:henkinset} Let $\sig = \langle \C, \mathcal{P}
	\rangle$ and $\Sigma \subseteq Sent(\sig)$. $\Sigma$ is an
	\sig-\textit{Henkin set} if and only if for every $B \in Form(\sig)$ and
	$x \in \mathcal{V}$, if $\Sigma \vdash \exists xB$, then $\Sigma \vdash
	B(c/x) \land \D c$, for some $c \in \C$.

	\end{df}

	\begin{df}\label{df:regularset} Let $\sig = \langle \C, \mathcal{P}
	\rangle$ and $\Sigma \subseteq Sent(\sig)$. $\Sigma$ is an
	\sig-\textit{regular set} if and only if (i) $\Sigma$ is non-trivial,
	(ii) $\Sigma$ is closed under $\vdash$, and (iii) $\Sigma$ is a
	disjunctive set, i.e., $\Sigma \vdash A \lor B$ only if $\Sigma \vdash
	A$ or $\Sigma \vdash B$.

	\end{df}

	\begin{lm}\label{lm:dprops} Let $\sig = \langle \C, \mathcal{P} \rangle$
	and $\Sigma \subseteq Sent(\sig)$. If $\Sigma$ is an \sig-regular Henkin
	set, then:

		\begin{itemize}\setl

		\item[1.] $B \land C \in \Sigma$ iff $B \in \Sigma$ and $C \in
		\Sigma$;

		\item[2.] $B \lor C \in \Sigma$ iff $B \in \Sigma$ or $C \in
		\Sigma$;

		\item[3.] $\neg (B \land C) \in \Sigma$ iff $\neg B \in \Sigma$
		or $\neg C \in \Sigma$;

		\item[4.] $\neg (B \lor C) \in \Sigma$ iff $\neg B \in \Sigma$
		and $\neg C \in \Sigma$;

		\item[5.] $\neg \neg B \in \Sigma$ iff $B \in \Sigma$;

		\item[6.] $\exists xB \in \Sigma$ iff $B(c/x), \D c \in \Sigma$, for
		some $c \in \C$;

		\item[7.] $\neg \forall xB \in \Sigma$ iff $\neg B(c/x), \D c
		\in \Sigma$, for some $c \in \C$.

		\end{itemize}

	\end{lm}

Lemma \ref{lm:extension} below is the first step of the proof of the
completeness of \ffde. Given a set of sentences $\Gamma$ that does not derive a
sentence $A$, it shows how to extend $\Gamma$ to a Henkin regular set $\Delta$
that still does not derive $A$.

	\begin{lm}\label{lm:extension} Let $\sig = \langle \C, \mathcal{P}
	\rangle$ and $\Gamma \cup \{A\} \subseteq Sent(\sig)$. If $\Gamma
	\nvdash A$, then there is a language $\sig^{+} = \langle \C^{+},
	\mathcal{P} \rangle$ and an $\sig^{+}$-regular Henkin set $\Delta$ such
	that $\C \subseteq \C^{+}$, $\Gamma \subseteq \Delta$, and $\Delta
	\nvdash A$.

	\end{lm}

	\begin{pr} Let $\C^{+} = \C \cup \{c_{i}: i \in \mathbb{N}\}$, with $\{c_{i}: i
	\in \mathbb{N}\} \cap (\C \cup \mathcal{P}) = \emptyset$, and $\sig^{+} =
	\langle \C^{+}, \mathcal{P} \rangle$. Adopt a fixed enumeration
	$B_{0},B_{1},B_{2},\dots$ of the sentences in $Sent(\sig^{+})$, and
	define the sequence $\langle j_{n} \rangle_{n \in \mathbb{N}}$ of natural
	numbers as follows: $j_{0} =$ the least natural number $k$ such that
	$c_{k}$ does not occur in $B_{0}$; and $j_{n+1} =$ the least natural
	number $k$ such that $c_{k}$ does not occur in $B_{n+1}$ and for every
	$i \leq n$, $k \neq j_{i}$.

	\vspace{0.25cm} \noi Now, consider the sequence $\langle \Gamma_{n}
	\rangle_{n \in \mathbb{N}}$ defined by:

		\begin{itemize}\setl

		\item $\Gamma_{0} = \Gamma$;

		\item $\Gamma_{n+1} = \left\{\begin{array}{lll} \Gamma_{n}	& \mbox{if } \Gamma_{n}, B_{n} \vdash A										\\

											&													\\

			\Gamma_{n} \cup \{B_{n}\} 					& \mbox{if } \Gamma_{n}, B_{n} \nvdash A \mbox{ and } B_{n} \neq \exists xC				\\
											& \mbox{ for every } C \in Form(\sig^{+})								\\
											&													\\
			\Gamma_{n} \cup \{B_{n}, C(c_{j_{n}}/x), \D c_{j_{n}} \} 	& \mbox{ if } \Gamma_{n}, B_{n} \nvdash A \mbox{ and } B_{n} = \exists xC				\\
											& \mbox{ for some } C \in Form(\sig^{+}) \mbox{ which has at most } x \mbox{ free}  			\\

		\end{array} \right.$

		\end{itemize}

	\noi Let $\Delta = \bigcup_{n \in \mathbb{N}}\Gamma_{n}$. Clearly, $\Gamma
	\subseteq \Delta$. We shall now prove that $\Delta$ is a regular Henkin
	set such that $\Delta \nvdash A$. It suffices to prove the following
	facts:

		\begin{itemize}\setl

		\item[(1)] For every $n \in \mathbb{N}$, $\Gamma_{n} \nvdash A$:
		The proof proceeds by induction on $n$.

		\item[(2)] $\Delta \nvdash A$: immediate consequence of (1).

		\item[(3)] If $\Delta \vdash C$, then $C \in \Delta$: immediate
		consequence of (2).

		\item[(4)] If $\Delta \vdash C \lor D$, then $\Delta \vdash C$ or
		$\Delta \vdash D$: immediate consequence of (2).\footnote{For
		detailed proofs of facts (1)-(4), see \cite[Lemma
		14]{Antunes&Rodrigues2024}.}

		\item[(5)] For every $C \in Form(\sig^{+})$, if $\Delta \vdash
		\exists xC$, then $\Delta \vdash C(c/x) \land \D c$, for some $c
		\in \C^{+}$: Let $n \in \mathbb{N}$ be such that $\exists xC =
		B_{n}$ and suppose that $\Delta \nvdash C(c/x) \land \D c$, for
		every $c \in \C^{+}$. In particular, $\Delta \nvdash
		C(c_{j_{n}}/x) \land \D c_{j_{n}}$. It then follows from the
		definition of $\langle \Gamma_{n} \rangle_{n \in \mathbb{N}}$
		that $\Gamma_{n}, B_{n} \vdash A$. Hence, $\Delta, B_{n} \vdash
		A$, and so $\Delta \nvdash B_{n}$ (i.e., $\Delta \nvdash \exists
		xC$), for otherwise $A$ would be derivable from $\Delta$.

		\end{itemize}

	\end{pr}

The following technical lemma will be necessary below:

	\begin{lm}\label{lm:newconstants} Let $\sig = \langle \C, \mathcal{P}
	\rangle$ and $\sig' = \langle \C', \mathcal{P} \rangle$ be such that $\C
	\subseteq \C'$. If $\Gamma \cup \{B\} \subseteq Sent(\sig)$, then
	$\Gamma \vdash_{\sig} B$ if and only if $\Gamma \vdash_{\sig'} B$.

	\end{lm}

Lemma \ref{lm:canonical} below, which corresponds to the second step of the
completeness proof, shows how to construct a model out of a regular Henkin set.
As usual, the model is defined in terms of items construed over the language and
the derivability relation.

	\begin{lm}\label{lm:canonical} Let $\sig = \langle \C, \mathcal{P}
	\rangle$ and $\Delta \cup \{A\} \subseteq Sent(\sig)$. If $\Delta$ is an
	\sig-regular Henkin set such that $\Delta \nvdash A$, then there is a
	model \K\ for \sig\ and a stage $w$ of \K\ such that $\K, w \vDash
	\Delta$ and $\K, w \nvDash A$.

	\end{lm}

	\begin{pr} Define the sequence $\langle \sig_{n} \rangle_{n \in
	\mathbb{N}}$ of first-order languages as follows:

		\begin{itemize}\setl

		\item $\sig_{0} = \sig$;

		\item $\sig_{n+1} = \langle \C_{n+1}, \mathcal{P} \rangle$,
		where $\C_{n+1} = \C_{n} \cup \{c^{n+1}_{i}: i \in \mathbb{N}\}$ and
		$\C_{n} \cap \{c^{n+1}_{i}: i \in \mathbb{N}\} = \emptyset$.

		\end{itemize}

	\noi Let $W$ be the set $\{\Sigma\!: \Sigma \supseteq \Delta \mbox{ and }
	\Sigma \mbox{ is an } \sig_{n}\mbox{-regular Henkin set, for some } n
	\in \mathbb{N}\}$ and let $\leq$ be the inclusion relation on $W$. Take
	$d$ to be the function on $W$ such that for every $n \in \mathbb{N}$ and
	every $\sig_{n}$-regular Henkin set $\Sigma$, $d(\Sigma) = \{c: c \in
	\C_{n} \mbox{ and } \D c \in \Sigma\}$. Let $i_{+}$ and $i_{-}$ be
	functions on $W$ such that $\langle c_{1},c_{2} \rangle \in
	i_{+}(\Sigma)$ iff $c_{1} \oid c_{2} \in \Sigma$ and $\langle
	c_{1},c_{2} \rangle \in i_{-}(\Sigma)$ iff $c_{1} \noid c_{2} \in
	\Sigma$, for every $c_{1},c_{2} \in d(\Sigma)$. Finally, let $\sig_{\K}
	= \langle \C_{\K}, \mathcal{P} \rangle$ be the language that is just
	like $\langle \bigcup_{n \in \mathbb{N}}\C_{n}, \mathcal{P} \rangle$
	except that if $c \in \C_{0} = \C$ is such that $\D c \in \bigcup W$,
	then, in addition to $c$, $\C_\K$ also includes a new constant $\ov{c}$,
	and every $c \in \bigcup_{n \geq 1}\C_{n}$ such that $\D c \in \bigcup
	W$ is replaced by $\ov{c}$.

	\vspace{0.25cm} \noi By letting $e, e_{1},e_{2},\dots$ range over
	$\C_{\K}$, we may define the partial function $\wh{I}$ as follows:

		\begin{itemize}\setl

		\item[(i)] If $\D e^{*} \in \Sigma$, then $\wh{I}(e, \Sigma) =
		e^{*}$, and undefined otherwise;

		\item[(ii)] For every $P \in \mathcal{P}$, $\langle
		\wh{I}(e_{1},\Sigma),\dots,\wh{I}(e_{m},\Sigma) \rangle \in
		P^{\Sigma}_{+}$ iff $Pe_{1}^{*}\dots e_{m}^{*} \in \Sigma$;

		\item[(iii)] For every $P \in \mathcal{P}$, $\langle
		\wh{I}(e_{1},\Sigma),\dots, \wh{I}(e_{m},\Sigma) \rangle \in
		P^{\Sigma}_{-}$ iff $\neg Pe_{1}^{*}\dots e_{m}^{*} \in \Sigma$;

		\item[(iv)] $\langle \wh{I}(e_{1},\Sigma),
		\wh{I}(e_{2},\Sigma) \rangle \in i_{+}(\Sigma)$ iff $e_{1}^{*}
		\oid e_{2}^{*} \in \Sigma$, and $\langle \wh{I}(e_{1},\Sigma),
		\wh{I}(e_{2},\Sigma) \rangle \in i_{-}(\Sigma)$ iff $e_{1}^{*}
		\noid e_{2}^{*} \in \Sigma$.

		\end{itemize}

	\noi where $e^{*} = c$ if $e = c \in \C_{0}$ or $e = \ov{c} \in
	\C_{\K}\setminus\C_{0}$.

	\vspace{0.25cm} \noi Consider next the structure $\wh{\K} = \langle W,
	\leq, d, \wh{I}, i_{+}, i_{-} \rangle$. Clearly, $\wh{I}$ satisfies all
	monotonicy conditions of Definition \ref{df:kripkemodel}(iv). By rules
	$\oid\!I$, $\oid\!E$, and $\noid\!I$, $i_{+}(\Sigma)$ and
	$i_{-}(\Sigma)$ satisfy conditions (a) and (b) of Definition
	\ref{df:kripkemodel}(v), and since $\leq$ is the subset relation on $W$,
	$i_{+}(\Sigma)$ and $i_{-}(\Sigma)$ also satisfy condition (c). Hence,
	$\wh{\K}$ is Kripke model for $\sig_{\K}$.

	\vspace{0.25cm} \noi Now, define the mapping $v: Sent(\sig_{\K}) \times
	W \longrightarrow \{0,1\}$ by $v(B,\Sigma) = 1 \mbox{ iff } B^{*} \in
	\Sigma$, where $B^{*}$ results from $B$ by replacing every occurrence of
	$e$ by $e^{*}$. We shall prove that $v$ satisfies all clauses of
	Definition \ref{df:valuation}. That $v$ satisfies clauses
	($v$7)-($v$11), ($v$13), and ($v$14) is an immediate consequence of
	Lemma \ref{lm:dprops}. Thus, we shall only prove that it satisfies
	clauses ($v$1), ($v$2), ($v$3), ($v$5), and ($v$12), and leave the
	proofs corresponding to ($v$4), ($v$6), and ($v$15) to the reader.

	\vspace{0.25cm} \noi Let $n \in \mathbb{N}$ and suppose that $\Sigma$ is
	an arbitrary $\sig_{n}$-regular Henkin set.

		\begin{itemize}\setl

		\item Clause ($v$1): Suppose that $v(\D e, \Sigma) = 1$. By the
		definition of $v$, $(\D e)^{*} = \D e^{*} \in \Sigma$. Hence,
		$\wh{I}(e, \Sigma)$ is defined. Suppose now that $\wh{I}(e,
		\Sigma)$ is undefined. Thus, $(\D e)^{*} \notin \Sigma$, and so
		$v(\D e, \Sigma) = 0$.

		\item Clause ($v$2): Suppose that $v(\neg \D e, \Sigma) = 1$.
		By the definition of $v$, $(\neg \D e)^{*} = \neg \D e^{*} \in
		\Sigma$. By rule $PEX_{\D}$ and the non-triviality of $\Sigma$,
		it then follows that $\D e^{*} \notin \Sigma$. Therefore,
		$\wh{I}(e, \Sigma)$ is undefined. Suppose now that $\wh{I}(e,
		\Sigma)$ is undefined. Thus, $\D e^{*} \notin \Sigma$, and so
		$\neg \D e^{*} \in \Sigma$, by rule $PEM_{\D}$ and Lemma
		\ref{lm:dprops}(2). Therefore, $v(\neg \D e, \Sigma) = 1$.

		\item Clause ($v$3): Suppose that $v(e_{1} \oid e_{2}, \Sigma) =
		1$. By the definition of $v$, $e_{1}^{*} \oid e_{2}^{*} \in
		\Sigma$, and so $\D e_{1}^{*}, \D e_{2}^{*} \in \Sigma$ (by rule
		$\D I$). Thus, $\wh{I}(e_{1}, \Sigma)$ and $\wh{I}(e_{2},
		\Sigma)$ are defined and $\langle \wh{I}(e_{1}, \Sigma),
		\wh{I}(e_{2}, \Sigma) \rangle = \langle e_{1}^{*}, e_{2}^{*}
		\rangle \in i_{+}(\Sigma)$. Suppose now that
		$\wh{I}(e_{1},\Sigma)$ and $\wh{I}(e_{2}, \Sigma)$ are defined
		and $\langle \wh{I}(e_{1}, \Sigma), \wh{I}(e_{2}, \Sigma)
		\rangle \in i_{+}(\Sigma)$. Hence, $\langle e_{1}^{*},
		e_{2}^{*} \rangle \in i_{+}(\Sigma)$, and so $e_{1}^{*} \oid
		e_{2}^{*} \in \Sigma$. Therefore, $v(e_{1} \oid e_{2}, \Sigma) =
		1$.

		\item Clause ($v$5): Note first that $Pe_{1}^{*}\dots e_{m}^{*}
		= (P e_{1}\dots e_{m})^{*}$. Suppose that $\wh{I}(e_{i},\Sigma)$
		is defined, for every $1 \leq i \leq m$, and that $\langle
		\wh{I}(e_{1},\Sigma),\dots,\wh{I}(e_{m},\Sigma) \rangle \in
		P^{\Sigma}_{+}$. Thus, $Pe_{1}^{*}\dots e_{m}^{*} = (Pe_{1}\dots
		e_{m})^{*} \in \Sigma$, and so $v(Pe_{1}\dots e_{m},\Sigma) =
		1$. Suppose, on the other hand, that $\wh{I}(e_{i},\Sigma)$ is
		undefined, for some $1 \leq i \leq m$. Thus, $\D e_{i}^{*}
		\notin \Sigma$, and so $(Pe_{1}\dots e_{m})^{*} \notin \Sigma$
		(by rule $\D I$). Therefore, $v(Pe_{1}\dots e_{m}, \Sigma) =
		0$. Finally, assume that $\wh{I}(e_{i},\Sigma)$ is defined, for
		every $1 \leq i \leq m$, but that $\langle
		\wh{I}(e_{1},\Sigma),\dots,\wh{I}(e_{m},\Sigma) \rangle \notin
		P^{\Sigma}_{+}$. Hence, $(Pe_{1}\dots e_{m})^{*} \notin \Sigma$,
		and so $v(Pe_{1}\dots e_{m},\Sigma) = 0$.

		\item Clause ($v$12): Note first that $\forall xB \in \Sigma$ is
		equivalent to the following claim:

			\begin{itemize}

			\item[($\dagger$)] For every $k \geq n$ and every
			$\sig_{k}$-regular Henkin set $\Sigma' \in W$ such that
			$\Sigma' \supseteq \Sigma$, $B(c/x) \in \Sigma'$, for
			every $c \in \C_{k}$ such that $\D c \in \Sigma'$.

			\end{itemize}

		For suppose that $\forall xB \in \Sigma$. Hence, $\forall xB \in
		\Sigma'$, since $\Sigma \subseteq \Sigma'$. Let $c \in \C_{k}$
		be such that $\D c \in \Sigma'$. It then follows by rule
		$\forall E_{\D}$ that $\Sigma' \vdash_{\sig_{k}} B(c/x)$.
		Therefore, $B(c/x) \in \Sigma'$, for every $c \in \C_{k}$ such
		that $\D c \in \Sigma'$ (since $c$ was arbitrary). As for the
		other direction, suppose that $\forall xB \notin \Sigma$ and let
		$c \in \C_{n+1}/\C_{n}$. Assume that $\Sigma, \D c
		\vdash_{\sig_{n+1}} B(c/x)$. Since $c$ does not occur in $\Sigma
		\cup \{B\}$, it follows by rule $\forall I_{\D}$ that $\Sigma
		\vdash_{\sig_{n+1}} \forall xB$. By Lemma \ref{lm:newconstants},
		$\Sigma \vdash_{\sig_{n}} \forall xB$, and so $\forall xB \in
		\Sigma$. But this result contradicts the initial hypothesis.
		Hence, $\Sigma, \D c \nvdash_{\sig_{n+1}} B(c/x)$. By taking
		$\Gamma$ to be $\Sigma \cup \{\D c\}$ and $\C^{+}$ to be
		$\C_{n+1} \cup \{c^{n+2}_{i}: i \in \mathbb{N}\}$ in the proof
		of Lemma \ref{lm:extension}, it then follows that there is an
		$\sig_{n+2}$-regular Henkin set $\Sigma' \supseteq \Sigma$ such
		that $B(c/x) \notin \Sigma'$, even though $\D c \in \Sigma'$.

		\vspace{0.25cm} \noi Now, suppose that $v(\forall xB, \Sigma) =
		1$. By the definition of $v$, $(\forall xB)^{*} = \forall xB^{*}
		\in \Sigma$. Hence, for every $k \geq n$ and every $\Sigma'
		\supseteq \Sigma$, $B^{*}(c/x) \in \Sigma'$, for every $c \in
		\C_{k}$ such that $\D c \in \Sigma'$ (by ($\dagger$)). Since
		$\ov{c}^{*} = c$, $(B(\ov{c}/x))^{*} = B^{*}(\ov{c}^{*}/x) =
		B^{*}(c/x)$. It then follows that for every $k \geq n$ and every
		$\Sigma' \supseteq \Sigma$, $v(B(\ov{c}/x), \Sigma') = 1$, for
		every $c \in d(\Sigma')$. Suppose, on the other hand, that
		$v(\forall xB, \Sigma) = 0$. By the definition of $v$, $(\forall
		xB)^{*} = \forall xB^{*} \notin \Sigma$. Hence, there is a
		natural number $k \geq n$ and an $\sig_{k}$-regular Henkin set
		$\Sigma' \supseteq \Sigma$ such that $B^{*}(c/x) \notin
		\Sigma'$, for some $c \in \C_{k}$ such that $\D c \in \Sigma'$
		(by ($\dagger$)). Therefore, there is $c \in d(\Sigma')$ such
		$v(B(\ov{c}/x), \Sigma') = 0$.

		\end{itemize}

	\noi Finally, let $I$ be the restriction of $\wh{I}$ to the parameters
	in $\sig = \sig_{0}$ and let $\K$ be the model $\langle W, \leq, d, I,
	i_{+}, i_{-} \rangle$. Clearly, \K\ is a model for \sig, $\sig_{\K}$ is
	its diagram language, $\wh{\K}$ is its associated structure, and $v$ is
	the valuation function induced by $\K$. Given the way $v$ was defined,
	it follows that $v(B, \Delta) = 1$ if and only $B^{*} \in \Delta$, for
	every $B \in Sent(\sig)$. But because none of the sentences of \sig\
	have occurrences of the $\ov{c}$ constants, it follows that $v(B,
	\Delta) = 1$ if and only if $B \in \Delta$, for every $B \in
	Sent(\sig)$. Therefore, there is $w \in W$ (viz., $\Delta$ itself) such
	that $\K, w \vDash \Delta$ and $\K, w \nvDash A$.

	\end{pr}

	\begin{tr}\label{tr:completeness} Let $\Gamma \cup \{A\} \subseteq
	Sent(\sig)$. If $\Gamma \vDash A$, then $\Gamma \vdash_{\sig} A$.

	\end{tr}

	\begin{pr} Suppose that $\Gamma \nvdash A$. By Lemma \ref{lm:extension},
	there is a language $\sig^{+} = \langle \C^{+}, \mathcal{P} \rangle$ and
	an $\sig^{+}$-regular Henkin set $\Delta$ such that $\C \subseteq
	\C^{+}$, $\Gamma \subseteq \Delta$, and $\Delta \nvdash A$. By Lemma
	\ref{lm:canonical}, there is a model $\K = \langle W, \leq, d, I, i_{+},
	i_{-} \rangle$ for $\sig^{+}$ and a stage $w \in W$ such that $\K, w
	\vDash \Delta$ and $\K, w \nvDash A$. Thus, $\K, w \vDash \Gamma$ and
	$\K, w \nvDash A$. Let $\K_{\sig} = \langle W, \leq, d, I_{\sig}, i_{+},
	i_{-} \rangle$ be such that $I_{\sig}$ is the restriction of $I$ to the
	parameters in $\sig$. Clearly, $\K_{\sig}$ is a model for $\sig$ such
	that $\K_{\sig}, w \vDash \Gamma$ and $\K_{\sig}, w \nvDash A$. Hence,
	$\Gamma \nvDash A$.

	\end{pr}

	\begin{rmrk}\label{rmrk:congruence} By taking a closer look at the proof
	of Lemma \ref{lm:canonical} we can see why the strategy adopted above
	for proving the completeness of \ffde\ would not work had the extension
	of $\oid$ at each stage $w$ been interpreted as the identity relation on
	$d(w)$. According to the standard strategy for proving the completeness
	of first-order systems with identity, the domain of each stage $\Sigma$
	would be defined as the set of equivalence classes $[c]_{\Sigma}$ such
	that $[c]_{\Sigma} = \{c'\!: c \oid c' \in \Sigma\}$. However, we would
	then be prevented from ensuring that the domains of stages are preserved
	across $\leq$-related stages, that is, that $d(\Sigma) \subseteq
	d(\Sigma')$ whenever $\Sigma \subseteq \Sigma'$. Suppose, for instance,
	that $c \in \C_{n}$ and that $\Sigma$ is an $\sig_{n}$-regular Henkin
	set. Suppose further that $c' \in \C_{n+1}$ and that $\Sigma'$ is an
	$\sig_{n+1}$-regular Henking set such that $\Sigma \subseteq \Sigma'$
	and $c \oid c' \in \Sigma'$. Since $c \oid c' \notin \Sigma$,
	$[c]_{\Sigma}$ is different from $[c]_{\Sigma'}$. Therefore,
	$[c]_{\Sigma}$ does not belong to the domain of $\Sigma'$, even though
	$\Sigma \leq \Sigma'$.

	\end{rmrk}

\section{The logic \fnf}\label{sec:fn4}

{\az deletar: As \fde, \ffde\ also lacks an implication that validates \textit{modus ponens}.
In order to overcome this limitation, in this section we extend \ffde\ with a
constructive implication $\rightarrow$, obtaining the logic \fnf. \fnf\ is a
universally free version with identity of the logic \qnf\ (a.k.a. \nm), \textbf{the
variable domain paraconsistent and paracomplete first-order (N4 é usualmente chamada paraconsistent nelson logic, fragmento de \textit{N}) 
logic introduced by
Almukdad and Nelson in \cite{Almukdad&Nelson1984}. The propositional fragment of
\fnf\ is usually referred to as \textit{N4}, which is the paraconsistent version
of Nelson's well-known logic of constructive negation \nt\ typo: N3 (parece q a prop é N3) \cite{Nelson1949}.}
In \cite{Antunes&Rodrigues2024} we showed how \qnf\ can be obtained from \qvfde\
and presented there a corresponding natural deduction system that is sound and
complete with respect to an appropriate Kripke semantics. We will now show how
\ffde\ can be modified in a similar vein to yield \fnf.
}

{\red As \fde, \ffde\ also lacks an implication that validates \textit{modus ponens}.
In order to overcome this limitation, in this section we extend \ffde\ with a
constructive implication $\rightarrow$, obtaining the logic \fnf. \fnf\ is a
universally free version with identity of the logic \qnf\ (a.k.a. \nm), the
paraconsistent 
first-order logic introduced by
Almukdad and Nelson in \cite{Almukdad&Nelson1984}.  
In \cite{Antunes&Rodrigues2024} we showed how \qnf\ can be obtained from \qvfde\
and presented there a corresponding natural deduction system that is sound and
complete with respect to an appropriate Kripke semantics. We will now show how
\ffde\ can be modified in a similar vein to yield \fnf.
}

 \begin{df}\label{df:fnel} \fnf\ is obtained by adding the following
	rules to \ffde:

		{\footnotesize \begin{center}

				\bottomAlignProof
							\AxiomC{$[A]$} \noLine
						\UnaryInfC{$\vdots$} \noLine
						\UnaryInfC{$B$}
					\RightLabel{$\rightarrow I$}
					\UnaryInfC{$A \rightarrow B$}
				\DisplayProof
			\qquad
				\bottomAlignProof
						\AxiomC{$A \rightarrow B$}
						\AxiomC{$A$}
					\RightLabel{$\rightarrow E$}
					\BinaryInfC{$B$}
				\DisplayProof
			\qquad
				\bottomAlignProof
						\AxiomC{$A$}
						\AxiomC{$\neg B$}
					\RightLabel{$\neg\!\rightarrow I$}
					\BinaryInfC{$\neg (A \rightarrow B)$}
				\DisplayProof
			\qquad
				\bottomAlignProof
						\AxiomC{$\neg(A \rightarrow B)$}
					\RightLabel{$\neg\!\rightarrow E$}
					\UnaryInfC{$A$}
				\DisplayProof
			$\!\!$
				\bottomAlignProof
						\AxiomC{$\neg(A \rightarrow B)$}
					\UnaryInfC{$\neg B$}
				\DisplayProof

		\end{center}

		}

	\end{df}

	\begin{df}\label{df:kripkemodel-n} Let $\sig = \langle \C, \mathcal{P}
	\rangle$ be a first-order language and $\K = \langle W, \leq, d, I,
	i_{+}, i_{-} \rangle$ be an \ffde-model for \sig. \K\ is also an
	\fnf-\textit{model}, but in addition to clauses ($v$1)-($v$15), its
	valuation function $v$ also satisfies:

		\begin{itemize}\setl

		\item[($v$16)] $v(A \rightarrow B, w) = 1$ iff for every $w' \geq
		w$, if $v(A, w') = 1$, then $v(B, w') = 1$;

		\item[($v$17)] $v(\neg(A \rightarrow B, w)) = 1$ iff $v(A, w) = 1$
		and $v(\neg B, w) = 1$.

		\end{itemize}

	\end{df}

	\begin{lm}\label{lm:persistence-n} Let $A \in Sent(\sig)$ and suppose
	that \K\ is an \fnf-model for \sig. For every stages $w$ and $w'$ of
	\K, if $\K, w \vDash A$ and $w \leq w'$, then $\K, w' \vDash A$.

	\end{lm}

	\begin{tr}\label{tr:completeness-n} \textbf{(Soundness and
	Completeness)} The natural deduction system for \fnf\ above is sound
	and complete with respect to the class of all \fnf-models.

	\end{tr}

	\begin{pr} \textbf{(Soundness)} We shall only care about rule
	$\rightarrow I$, since proving the soundness of rules $\rightarrow E$,
	$\neg\!\rightarrow I$, and $\neg\!\rightarrow E$ is straightforward.
	Assume that $\Gamma, A \vDash B$ and that $\K, w \vDash \Gamma$. Let $w'
	\geq w$ and suppose that $v(A, w') = 1$. By Lemma
	\ref{lm:persistence-n}, it then follows that $\K, w' \vDash \Gamma \cup
	\{A\}$. Therefore, $v(B, w') = 1$. This result, together with the
	soundness proof for \ffde, suffices for establishing the soundness of
	\fnf. \\

	\noi \textbf{(Completeness)} Note first that if $\Sigma$ is a regular
	Henkin set, then:

		\begin{itemize}\setl

		\item[1.] $B \rightarrow C \in \Sigma$ iff for every regular
		Henkin set $\Sigma' \supseteq \Sigma$, $B \in \Sigma'$ only if
		$C \in \Sigma$; and

		\item[2.] $\neg (B \rightarrow C) \in \Sigma$ iff $B \in \Sigma$
		and $\neg C \in \Sigma$.

		\end{itemize}

	\noi For if $B \rightarrow C \in \Sigma$ and $B \in \Sigma'$, then $B
	\rightarrow C, B \in \Sigma'$. Thus, $\Sigma' \vdash C$, and so $C \in
	\Sigma'$. If, on the other hand, $B \rightarrow C \notin \Sigma$, then
	$\Sigma, B \nvdash C$ (by rule $\rightarrow I$). By (the \fnf-analogue
	of) Lemma \ref{lm:extension}, whose proof is exactly the same as before,
	there is a regular Henkin set $\Sigma' \supseteq \Sigma \cup \{B\}$ such
	that $\Sigma' \nvdash C$. Thus, $B \in \Sigma'$ and $C \notin
	\Sigma'$.\footnote{Item (2) above is an immediate consequence of rules
	$\neg \! \rightarrow I$ and $\neg\!\rightarrow E$ and the properties of
	regular Henkin sets.}

	\vspace{0.25cm} \noi Facts (1) and (2) above ensure that the function
	$v$ defined in the proof of Lemma \ref{lm:canonical} would satisfy
	clauses ($v$16) and ($v$17) in the proof of the corresponding result for
	\fnf, and so that $v$ would also count as a valuation in the relevant
	sense. These observations suffice to make it clear how the completeness
	proof for \ffde\ presented above can be easily adapted to the case of
	\fnf.

	\end{pr}

\section{Final remarks}\label{sec:final}

The aim of this paper was to introduce quantified versions of Belnap-Dunn's
four-valued logic, or the logic of first-degree entailment, and Nelson's
paraconsistent logic \nf, remaining as faithful as possible to Belnap-Dunn's
proposal of an information-based logic. The logics presented here, \ffde\ and
\fnf, are variable domain universally free versions with identity of \qvfde\ and
\qnf. Like \qvfde\ and \qnf, \ffde\ and \fnf\ are capable of representing
positive and negative information in such a way that they are independent from
each other. This behavior manifests itself with respect to all logical symbols,
including the identity predicate. Moreover, as \qvfde\ and \qnf, \ffde\ and
\fnf\ are also capable of representing the dynamic aspect of information
processing, i.e., how a database receives new information over time.  However,
in addition to including the identity predicate, they improve upon their
non-free versions in at least two respects: by allowing for models that contain
empty stages, \ffde\ and \fnf\ are able to represent scenarios in which a
database acknowledges no individual whatsoever, a situation that may take place
when the database has not yet been fed up with any information; and by allowing
for empty constants, they may represent scenarios in which the linguistic
resources available exceed the individuals it acknowledges at a certain stage of
its development, such as when there is a predetermined way of generating names
that may or may not be assigned to individuals (e.g., social security numbers).

It should be observed at this point that in spite of the advantages mentioned
above \ffde\ and \fnf\ don't go far enough, for they are unable to countenance
the possibility of information being withdrawn from a database. \ffde\ and \fnf\
are Tarskian (and so monotonic) logics and their semantics satisfy the
persistence condition stated in Proposition \ref{prop:persistence} -- which,
when framed in terms of information, states that if a piece of information has
been stored in a database at a certain stage it cannot be detracted at further
stages. This limitation could be overcome by theories of belief revision based
on \ffde\ (\fnf), but devising such theories is a task we will leave for another
occasion (and perhaps for other researchers).

 {\red   \textbf{\az Nao sei se vc desconsiderou ou nao viu essa alteracao nos final remarks}  Before closing the paper, two comments are in order. 
First, it should be noted that this paper does not aim to propose or discuss concrete applications of the logics under investigation, e.g. in database management systems. 
As we have seen, Belnap's  motivation for introducing the four-valued semantics for \fde\ 
was to design an information-based logic, in the sense explained in Section \ref{sec:fde} above. 
 The logics \textit{FFDE} and \textit{FN4} are intended to be further developments of Belnap's  (and Dunn's) 
 ideas, in which configurations of a database (in Belnap's words, the computer's `epistemic states' 
 \cite[p.~39]{Belnap1977}) are thought of as possibly inconsistent and incomplete information states, 
 represented by worlds in Kripke models (as remarked, 
 this idea  traces back to Wansing \cite[pp.~15-16]{Wansing1993}). Our point here is that 
  universally free logics provide  an improvement on this interpretation in terms of information states.

Second, despite the advantages mentioned above, \ffde\ and \fnf\ fall short in one significant aspect: 
they do not allow for the possibility of withdrawing information from a database. 
These logics are Tarskian (and thus monotonic), meaning that their semantics satisfy the persistence condition stated in Proposition \ref{prop:persistence}. 
In terms of information, this condition implies that once a piece of information has been stored in a database at a certain stage, it cannot be retracted at later stages. 
This limitation could be overcome by theories of belief revision based
on \ffde\ and \fnf, but devising such theories is a task we will leave for
another occasion (and perhaps for other researchers). 
}

As we opened the paper with a passage from Belnap, we will close it with the
following passage from an interview given by Dunn in 2018:

	\begin{quote}

	[T]he main motivation that I see for paraconsistent logic comes from the
	Internet and its World Wide Web. Belnap was prescient in motivating what
	has become known as the Belnap-Dunn 4-valued logic by talking of an
	unstructured database in which various people might enter inconsistent
	data. This is before the World Wide Web. And now not only do we have
	this enormous unstructured database, but we have ``bots'' roaming it for
	information. Assuming that these bots will perform inferences, I see
	paraconsistent logic as an essential tool to prevent ``Explosion''.
	Ideally this would be combined with other tools for reasoning and
	learning such as probability and statistics, deep learning and other
	machine learning techniques\dots Perhaps the most theoretically
	challenging problem for paraconsistent logics has to do with
	incorporating belief revision and defeasible reasoning.
	\citep{Dunn2018}

	\end{quote}

\noi Dunn's words are an update of the ideas in \cite{Belnap1977} and mention
some of the developments that Benalp's original proposal must undergo,
incorporating belief revision and defeasible reasoning being among them. By
presenting variable domain universally free first-order extensions with identity
of \fde\ and \nf, we hope to have contributed to the task of improving Belnap's
proposal in other respects.


\bibliographystyle{plainnat}	
\bibliography{ffde}	

\end{document}